\documentclass{article}
\usepackage{latexsym}
\usepackage{epsfig}
\usepackage{amsmath}
\usepackage{amssymb}
\usepackage{amsthm}
\usepackage{caption}
\newtheorem{theorem}{Theorem}[section]
\newtheorem{corollary}[theorem]{Corollary}

\newtheorem{prop}[theorem]{Proposition}
\newtheorem{lemma}[theorem]{Lemma}

\newtheorem{fact}[theorem]{Fact}

\theoremstyle{remark}

\theoremstyle{definition}
\newtheorem{definition}[theorem]{Definition}

\DeclareMathOperator\Aut{Aut}

\DeclareMathOperator\diag{diag}

\DeclareMathOperator\tr{tr}

\begin{document}

\title{Self-associated three-dimensional cones}

\author{Roland Hildebrand \thanks{%
Univ.~Grenoble Alpes, CNRS, Grenoble INP, LJK, 38000 Grenoble, France
({\tt roland.hildebrand@univ-grenoble-alpes.fr}).}}

\maketitle

\begin{abstract}
For every proper convex cone $K \subset \mathbb R^3$ there exists a unique complete hyperbolic affine 2-sphere with mean curvature $-1$ which is asymptotic to the boundary of the cone. Two cones are associated if the corresponding affine spheres can be mapped to each other by an orientation-preserving isometry. This equivalence relation is generated by the groups $SL(3,\mathbb R)$ and $S^1$, where the former acts by linear transformations of the ambient space, and the latter by multiplication of the cubic holomorphic differential of the affine sphere by unimodular complex constants. The action of $S^1$ generalizes conic duality, which acts by multiplication of the cubic differential by $-1$. We call a cone self-associated if it is linearly isomorphic to all its associated cones, in which case the action of $S^1$ induces (nonlinear) isometries of the corresponding affine sphere. We give a complete classification of the self-associated cones and compute isothermal parametrizations of the corresponding affine spheres. The solutions can be expressed in terms of degenerate Painlev\'e III transcendents. The boundaries of generic self-associated cones can be represented as conic hulls of vector-valued solutions of a certain third-order linear ordinary differential equation with periodic coefficients, but there exist also cones with polyhedral boundary parts.
\end{abstract}

\section{Introduction}

In this work we classify the self-associated convex cones in $\mathbb R^3$, by which we mean those cones that are linearly isomorphic to all its associated cones. The notion of associated cones has been introduced in the paper \cite{LinWang16} of Z.~Lin and E.~Wang and is by virtue of the Calabi theorem derived from the notion of associated families of affine spheres \cite{SimonWang93}. The results in this paper in spirit partly resemble the work \cite{DumasWolf15} of Dumas and Wolf, which were the first to establish explicit non-trivial relations between affine spheres and their asymptotic cones.

Self-associated cones are akin to semi-homogeneous cones, which have been introduced in \cite{Hildebrand14b} and whose corresponding affine spheres' isothermal parametrizations have been explicitly computed in \cite{LinWang16}. While the affine spheres of the latter possess a continuous isometry group which is generated by linear automorphisms of the cone, the affine spheres of the former have isometries which are of a nonlinear nature and correspond to a continuous symmetry of the cone generalizing the duality symmetry of self-dual cones.


\subsection{Background}

A proper, or regular, convex cone $K \subset \mathbb R^n$ is a closed convex cone with non-empty interior and containing no line. In the sequel we shall speak of cones for brevity, meaning always proper convex cones. The \emph{Calabi conjecture} on affine spheres \cite{Calabi72}, proven by the efforts of several authors, states that for every cone there exists a unique complete hyperbolic affine sphere with mean curvature $-1$ which is asymptotic to the boundary of the cone \cite{ChengYau77,Sasaki80}, and conversely, every complete hyperbolic affine sphere is asymptotic to the boundary of a cone \cite{ChengYau86,Li90,Li92}. These affine spheres are equipped with a complete non-positively curved Riemannian metric, the affine metric $h$, and a totally symmetric trace-less $(0,3)$-tensor, the cubic form $C$. The two objects are equivariant under the action of the group $SL(n,\mathbb R)$ of unimodular transformations. Moreover, the affine sphere and the corresponding cone $K$ can be reconstructed from $(h,C)$ up to a unimodular transformation of the ambient space $\mathbb R^n$. For more on affine spheres see, e.g., \cite{NomizuSasaki}, for a survey on the Calabi conjecture see \cite[Chapter 2]{LiSimonZhao}.

Complete hyperbolic affine 2-spheres are non-compact simply connected Riemann surfaces, and their affine metrics have the form $h = e^u|dz|^2$ in an isothermal complex coordinate $z = x + iy \in M$, $M \subset \mathbb C$ being a simply connected domain and $u: M \to \mathbb R$ an analytic conformal factor. The cubic form can be represented as $C = 2Re(U(z)dz^3)$, with $U: M \to \mathbb C$ a holomorphic function, the so-called cubic differential (more precisely, the cubic differential is $U\,dz^3$, but we shall refer to $U$ for brevity). It satisfies the compatibility condition \cite{Wang91}
\begin{equation} \label{Wangs_equation}
e^u = \frac12\Delta u + 2|U|^2e^{-2u},
\end{equation}
also called \emph{Wang's equation}. Here $\Delta u = u_{xx} + u_{yy} = 4u_{z\bar z}$ is the Laplacian of $u$. In the sequel, when we speak of complete solutions $(u,U)$ of \eqref{Wangs_equation}, we always assume that the metric $h$ is complete and $U$ is holomorphic.

The affine sphere is given by an embedding $f: M \to \mathbb R^3$, which is determined by the solution $(u,U)$ up to a unimodular transformation of the ambient space $\mathbb R^3$. We shall speak of the complete solution $(u,U)$ of Wang's equation \eqref{Wangs_equation} as \emph{corresponding} to the affine sphere $f$ or the convex cone $K$. In \cite{Wang91} C.~Wang deduced the moving frame equations whose integration allows to reconstruct the affine sphere $f$ from a given solution $(u,U)$. 
Define the $3 \times 3$ real matrix $F = (e^{-u/2}f_x,e^{-u/2}f_y,f)$. The first two columns of $F$ form an $h$-orthonormal basis of the tangent space to the embedding, while the third column is the position vector of the embedding. The condition that the mean curvature of the affine sphere equals $-1$ is equivalent to unimodularity of $F$. 
From the structure equations \cite{SimonWang93}
\[ f_{zz} = u_zf_z - Ue^{-u}f_{\bar z},\quad f_{z\bar z} = \frac12e^uf,\quad f_{\bar z\bar z} = -\bar Ue^{-u}f_z + u_{\bar z}f_{\bar z}
\]
we obtain the frame equations
\begin{equation} \label{frame_equations}
\begin{aligned}
F_x &= 
F\begin{pmatrix} -e^{-u}Re\,U & \frac{u_y}{2} + e^{-u}Im\,U & e^{u/2} \\ -\frac{u_y}{2} + e^{-u}Im\,U & e^{-u}Re\,U & 0 \\ e^{u/2} & 0 & 0 \end{pmatrix},\\
F_y &= 
F\begin{pmatrix} e^{-u}Im\,U & -\frac{u_x}{2} + e^{-u}Re\,U & 0 \\ \frac{u_x}{2} + e^{-u}Re\,U & -e^{-u}Im\,U & e^{u/2} \\ 0 & e^{u/2} & 0 \end{pmatrix}.
\end{aligned}
\end{equation}
Choosing an arbitrary point $z_0 = x_0 + iy_0 \in M$ and an arbitrary initial value $F(z_0) = F_0 \in SL(3,\mathbb R)$, we can recover the embedding $f$ from the third column of the matrix-valued function $F(z)$ by integrating \eqref{frame_equations}.

The isothermal coordinate system is not unique, but defined up to conformal isomorphisms of the domain $M$. Let $\tilde M \subset \mathbb C$ be the pre-image of the domain $M$ under the biholomorphic map $b: w \mapsto z$, and consider the embedding $f \circ b: \tilde M \to \mathbb R^3$. Then the embeddings $f \circ b$ and $f$ define the same affine sphere, considered as a surface in $\mathbb R^3$, but the parametrization is different. Let $(\tilde u,\tilde U)$ be the solution of \eqref{Wangs_equation} corresponding to $f \circ b$. 
Invariance of the pair $(h,C)$ then yields the transformation law
\begin{equation} \label{Wang_transform}
\tilde U(w) = U(z)\left(\frac{dz}{dw}\right)^3,\qquad \tilde u(w) = u(z) + 2\log\left|\frac{dz}{dw}\right|.
\end{equation}

\medskip

Let us single out the following consequence for further reference.

{\fact \label{fact:basic_equivalence} Equivalence classes of complete solutions $(u,U)$ of \eqref{Wangs_equation} under bi-holomorphisms between the domains of definition are in one-to-one correspondence to equivalence classes of regular convex cones $K \subset \mathbb R^3$ under the action of the group $SL(3,\mathbb R)$. }

\medskip

Several questions arise:

1. Given a complete solution $(u,U)$ of \eqref{Wangs_equation}, describe the corresponding cone $K \subset \mathbb R^3$ and vice versa.

2. Given a holomorphic function $U$ on a simply connected non-compact domain $D \subset \mathbb C$, characterize all corresponding complete solutions $(u,U)$ of \eqref{Wangs_equation}.

3. Given an analytic function $u$ on a simply connected non-compact domain $D \subset \mathbb C$, characterize all corresponding complete solutions $(u,U)$ of \eqref{Wangs_equation}.

\medskip

Results on the first question are scarce. Dumas and Wolf \cite{DumasWolf15} have shown that if $U$ is a polynomial on $\mathbb C$, then the complete solution $(u,U)$ corresponds to a polyhedral cone in $\mathbb R^3$, with the degree of the polynomial being equal to the number of extreme rays of the cone less 3. On the other hand, a polyhedral cone $K$ corresponds to a solution $(u,U)$ with $U$ a polynomial on $\mathbb C$. 
In \cite{LinWang16} the solutions $(u,U)$ have been computed for the semi-homogeneous cones, i.e., for those cones which have a non-trivial continuous group of linear automorphisms.

\medskip

The second question has been answered completely. Obviously it suffices to consider the unit disc $\mathbb D$ and the complex plane $\mathbb C$ only, as any other simply connected non-compact domain is conformally equivalent to one of these two. Using techniques from \cite{WanAu94} Q.~Li observed that if $U$ is a cubic holomorphic differential on the unit disc, then \eqref{Wangs_equation} has a unique complete solution \cite{QLi19}. He used this to prove that if $U$ is a non-zero holomorphic cubic differential on $\mathbb C$, then \eqref{Wangs_equation} has a unique complete solution \cite[Theorem 3.1]{QLi19}. We may summarize these results as follows.

{\fact \label{fact:uniqueness} For every holomorphic $U$ on a simply-connected domain $M \subset \mathbb C$ except for the zero function on $M = \mathbb C$ there exists a unique complete solution $(u,U)$ of \eqref{Wangs_equation}. }

Prior to this, existence and uniqueness results have been obtained for cubic differentials $U$ induced on the universal covers of compact Riemann surfaces \cite{Wang91,Loftin01,Labourie07}, on the universal covers of punctured compact Riemann surfaces with poles at the punctures \cite{Loftin04}, for functions on $\mathbb D$ satisfying certain bounded-ness conditions \cite{BenoistHulin13,BenoistHulin14}, and for polynomials on $\mathbb C$ \cite{DumasWolf15}.

\medskip

The third question has been answered in \cite{SimonWang93}. 
A necessary and sufficient condition on the affine metric $h$ induced by the analytic function $u$ has been given to admit a solution $(u,U)$.
For any two such solutions the holomorphic functions $U$ differ by a multiplicative unimodular complex constant $e^{i\varphi}$, and multiplying the function $U$ of some solution by $e^{i\varphi}$ yields another solution. 

\medskip

The last result gives rise to the concept of associated cones which is central for the present paper. If an affine 2-sphere with mean curvature $-1$, affine metric $h = e^u|dz|^2$, and cubic differential $U$ is given, then the affine spheres constructed from the pairs $(u,e^{i\varphi}U)$, where $\varphi$ runs through $[0,2\pi)$, exhaust all affine 2-spheres with affine metric $h$ and mean curvature $-1$. The orbits of complete hyperbolic affine 2-spheres with respect to the action of $SL(3,\mathbb R)$ are hence arranged in 1-parametric \emph{associated families}, on which the circle group $S^1 \simeq \{ c \in \mathbb C \,|\, |c| = 1 \}$ acts by multiplication of the cubic differential $U$ by the unimodular group element. Affine spheres belonging to orbits in the same family are called \emph{associated} \cite{SimonWang93}.

By virtue of the Calabi theorem this notion can naturally be extended to cones in $\mathbb R^3$. The action of the circle group $S^1$ on the solutions of \eqref{Wangs_equation} induces an action on the set of $SL(3,\mathbb R)$-orbits of cones, and these $SL(3,\mathbb R)$-orbits are also arranged in 1-parametric families. 
Cones belonging to orbits in the same family are called \emph{associated} \cite{LinWang16}.


The action of the group $S^1$ on a single associated family of $SL(3,\mathbb R)$-orbits of cones does not need to be faithful. It may well be that two solutions $(u,U)$, $(u,cU)$ of Wang's equation \eqref{Wangs_equation} for $c \not= 1$ lead to isomorphic affine spheres or, equivalently, to the same $SL(3,\mathbb R)$-orbit of cones. In this contribution we characterize those cones $K \subset \mathbb R^3$ whose $SL(3,\mathbb R)$-orbit is a fixed point of the action of $S^1$, and compute the corresponding solutions $(u,U)$ of Wang's equation \eqref{Wangs_equation}. 

{\definition \label{def:self_asso} Let $K \subset \mathbb R^3$ be a cone and let $(u,U)$ be a solution of \eqref{Wangs_equation} corresponding to $K$. For $\varphi \in [0,2\pi)$, let $K_{\varphi} \subset \mathbb R^3$ be a cone corresponding to the solution $(u,e^{i\varphi}U)$ of \eqref{Wangs_equation}. We call the cone $K$ \emph{self-associated} if for every $\varphi \in [0,2\pi)$ the cone $K_{\varphi}$ is linearly isomorphic to $K$. }

\medskip

Note that we did not restrict the linear isomorphisms between $K$ and $K_{\varphi}$ to be unimodular, and hence the condition in Definition \ref{def:self_asso} is a priori weaker than the condition that the $SL(3,\mathbb R)$-orbit of $K$ is a fixed point of the action of $S^1$. Later we shall show that these two conditions are actually equivalent.



Explicit results on associated families of cones are scarce. In \cite{LinWang16} the associated families have been computed for the semi-homogeneous cones. From the results of Dumas and Wolf \cite{DumasWolf15} it follows that any cone associated to a polyhedral cone is also polyhedral with the same number of extreme rays, and that the self-associated polyhedral cones are exactly the cones over the regular $n$-gons.

In \cite[Corollary 4.0.4]{LoftinThesis} Loftin observed that multiplying the cubic differential of an affine 2-sphere by $-1$ leads to the projectively dual affine 2-sphere. 
However, if an affine 2-sphere is asymptotic to the boundary of a convex cone $K \subset \mathbb R^3$, then the projectively dual affine 2-sphere is asymptotic to the boundary of the dual cone $K^* = \{ y \in \mathbb R^3 \,|\, x^Ty \geq 0\ \forall\ x \in K \}$. Hence $K^*$ is always associated to $K$, and its $SL(3,\mathbb R)$-orbit can be obtained from the $SL(3,\mathbb R)$-orbit of $K$ by the action of the group element $e^{i\pi} \in S^1$. 

Any self-associated cone $K$ must hence be self-dual, in the sense that it is linearly isomorphic to its dual $K^*$. The property of being self-associated is therefore stronger than self-duality. The subject of the paper is the description and classification of the self-associated cones and the corresponding solutions of \eqref{Wangs_equation}.

\subsection{Outline} \label{subs:outline}

In this work we make a step towards a better understanding of associated families of 3-dimensional cones. 
We provide a full classification of self-associated cones and explicitly describe their boundaries as well as the affine spheres which are asymptotic to these boundaries. We now outline the contents of the paper, sketch the strategy of the proofs, and summarize the results.

\medskip

In Section \ref{sec:Killing} we express the condition that a cone $K \subset \mathbb R^3$ is self-associated equivalently by the ana\-ly\-tic conditions \eqref{KillingConditionU},\eqref{KillingConditionu} on the corresponding solution $(u,U)$ of Wang's equation. These conditions state the existence of a Killing vector field $\psi$ on $M$ whose flow of conformal automorphisms of the domain $M$ multiplies the cubic differential by unimodular complex constants. Actually, the Killing condition \eqref{KillingConditionu} follows already from \eqref{KillingConditionU} and the weaker condition that $\psi$ generates a 1-parameter group of conformal automorphisms of the domain $M$.

\medskip

In Section \ref{sec:cubic_class} we classify the solutions $(\psi,U)$ of \eqref{KillingConditionU}, where $\psi$ generates a subgroup of automorphisms of $M$, up to conformal isomorphisms. We distinguish two cases. In the first case the domain can be transformed to an open disc with radius $R \in (0,+\infty]$, the generated conformal automorphisms are the rotations of the disc, and $U(z) = z^k$, $k \in \mathbb N$. In the second case the domain can be transformed to a vertical strip $(a,b) + i\mathbb R$, where $-\infty \leq a < b \leq +\infty$, the automorphisms are the vertical translations, and $U(z) = e^z$. We shall refer to these cases as the rotational and the translational case, respectively.

The trivial solution $U \equiv 0$ of \eqref{KillingConditionU} which exists for every $\psi$ corresponds to the Lorentz cone $L_3$ in the case of a hyperbolic domain.

\medskip

In Section \ref{sec:solving_Wang} we find for every pair $(\psi,U)$ the unique complete solution $u$ of Wang's equation. Its invariance with respect to the automorphisms generated by $\psi$ allows to reduce \eqref{Wangs_equation} to the degenerate Painlev\'e III equation \eqref{PainleveIII}. For each $(\psi,U)$ we characterize the corresponding solution $v(t)$ of \eqref{PainleveIII}, in particular, its asymptotics at the boundary of the interval of definition.

Other Painlev\'e III equations play a role in the description of constant or harmonic inverse mean curvature surfaces with radial symmetry \cite{BobenkoIts95},\cite{BobenkoEitnerKitaev97}, Amsler surfaces \cite{BobenkoKitaev98}, and in Smyth surfaces and indefinite affine spheres with intersecting straight lines \cite{BobenkoEitner00}.



\medskip

In Section \ref{sec:integr_frame} we integrate the frame equations and compute pieces of the boundary of the self-associated cones. In order to capture the asymptotics of the moving frame at the boundary of the domain we shall employ the technique of \emph{osculation maps} which has been introduced in \cite{DumasWolf15}. We compare the moving frame $F$ with an explicit diverging unimodular matrix $V$ such that the ratio $G = FV^{-1}$ has a  finite limit as the argument tends to the boundary. The matrix $V$ is chosen such that the last column of $G$ and a suitably scaled multiple of the position vector $f$ tend to the same non-zero point in the considered limit. This point must be a boundary point of the cone.

The moving frame equations on $F$ translate to similar partial differential equations on $G$ and in general reduce to a linear third-order ordinary differential equation (ODE) with periodic coefficients on the limit point as a function of the line of convergence. In these cases the corresponding boundary piece of $K$ can be described as the conic hull of a vector-valued solution of the ODE.

In the remaining cases the limit point does not continuously depend on the line along which the limit is taken, but rather yields a discrete sequence of points in $\mathbb R^3$. In these cases the boundary piece consists of planar faces. Exactly this situation has been encountered by Dumas and Wolf in \cite{DumasWolf15}.

\medskip

In Section \ref{sec:assemble} we assemble the boundary pieces to obtain the whole boundary of the cone. The main tool to accomplish this task are the linear automorphisms $\Sigma,T$ of the cone that are generated by complex conjugation of the domain and by its rotation by an angle $\frac{2\pi}{k+3}$ or translation by $2\pi i$, respectively. These symmetries generate a subgroup of automorphisms which is isomorphic to a finite dihedral group in the rotational case (Lemma \ref{lem:symR}) and to the infinite dihedral group in the translational case (Lemma \ref{lem:symT}).

According to the whether the trace of $T$ is smaller, equal, or greater than 3 we distinguish cones of elliptic, parabolic, and hyperbolic type (Definition \ref{def:types}). The first arises in the rotational case, the second in the translational case with $a = -\infty$, and the last in the translational case with $a$ finite.

\medskip

In Section \ref{sec:results} we summarize our classification in Theorems \ref{thm:ell}, \ref{thm:par}, and \ref{thm:hyp}, describing the three types of cones, respectively.

\medskip

In Section \ref{sec:conclusion} we pose some open questions and suggest directions for further research.

\section{Killing vectors and self-associated cones} \label{sec:Killing}

In this section we characterize affine spheres which are asymptotic to the boundary of self-associated cones by the existence of a 1-parametric subgroup of automorphisms of the domain which non-trivially multiplies the cubic differential by unimodular complex constants. The infinitesimal generator of the subgroup turns out to be a Killing vector field (i.e., generator of isometries) for the affine metric. 

This will be accomplished in two steps. First we show that if all associated cones $K_{\varphi}$ of a cone $K \subset \mathbb R^3$ lie in the same $GL(3,\mathbb R)$-orbit as $K$, then they must also lie in the same $SL(3,\mathbb R)$-orbit. In a second step we show the existence of the Killing vector field.

\subsection{Isomorphisms with negative determinant} \label{subs:negdet}

Suppose two cones $K,\tilde K$ are linearly isomorphic, but not in the same $SL(3,\mathbb R)$-orbit. Then there exists a linear map with determinant $-1$ taking $K$ to $\tilde K$. The next result shows how the corresponding solutions of \eqref{Wangs_equation} are related to each other.


{\lemma \label{lem:detM1} Let $K \subset \mathbb R^3$ be a cone and let $A \in GL(3,\mathbb R)$ be a linear map with determinant $-1$. Let $(u,U)$ be the solution of \eqref{Wangs_equation} corresponding to $K$, defined on a simply connected domain $M \subset \mathbb C$. Then a solution $(\tilde u,\tilde U)$ of \eqref{Wangs_equation} corresponding to the cone $A[K]$ is given by $\tilde u(z) = u(\bar z)$ and $\tilde U(z) = \overline{U(\bar z)}$ on the complex conjugate domain $\bar M$. }

\begin{proof}
The domain $\bar M$ is simply connected, $\tilde U$ is holomorphic, and the pair $(\tilde u,\tilde U)$ satisfies Wang's equation \eqref{Wangs_equation}. Moreover, the metric $\tilde h = e^{\tilde u}|dz|^2$ is complete on $\bar M$, because complex conjugation is an isometry between $M$ and $\bar M$. Hence $(\tilde u,\tilde U)$ is a solution of \eqref{Wangs_equation}.

Let $F: M \to SL(3,\mathbb R)$ be a solution of the frame equations \eqref{frame_equations} such that the embedding $f: M \to \mathbb R^3$ given by the third column of $F$ is asymptotic to $\partial K$. It is not hard to see that equations \eqref{frame_equations} are invariant under the substitution $(x,y,u,U,f) \mapsto (x,-y,u,\bar U,-f)$. Note that $-A \in SL(3,\mathbb R)$. The embedding $\tilde f: \bar M \to \mathbb R^3$ defined by $\tilde f(z) = (-A)(-f(\bar z)) = Af(\bar z)$ is then asymptotic to $\partial A[K]$ and gives rise to a solution $\tilde F: \bar M \to SL(3,\mathbb R)$ of the frame equations corresponding to $(\tilde u,\tilde U)$. Therefore the solution $(\tilde u,\tilde U)$ corresponds to the cone $A[K]$.
\end{proof}

We now show that for a self-associated cone, the isomorphisms between the associated cones can be assumed to have positive determinant.

{\lemma \label{lem:SLorbit} A cone $K \subset \mathbb R^3$ is self-associated if and only if all cones which are associated to $K$ are in the $SL(3,\mathbb R)$-orbit of $K$. }

\begin{proof}
The "if" direction follows from the definition of self-associated cones. Let us prove the "only if" direction.

Let $K \subset \mathbb R^3$ be self-associated. Then every associated cone $K_{\varphi}$ is in the $SL(3,\mathbb R)$-orbit of $K$ or in the $SL(3,\mathbb R)$-orbit of $-K$. If $K,-K$ are in the same $SL(3,\mathbb R)$-orbit, then the claim of the lemma follows. Let us therefore suppose that the $SL(3,\mathbb R)$-orbits of the cones $K,-K$ are distinct. Define the disjoint complementary subsets $S_+,S_- \subset S^1$ such that $e^{i\varphi} \in S_{\pm}$ if $K_{\varphi}$ is in the $SL(3,\mathbb R)$-orbit of $\pm K$, respectively.

From \eqref{Wang_transform} it follows by multiplication of the first equation with $e^{i\varphi}$ that if two solutions $(u,U)$, $(\tilde u,\tilde U)$ of Wang's equation are conformally equivalent, then for all $\varphi$ the solutions $(u,e^{i\varphi}U)$, $(\tilde u,e^{i\varphi}\tilde U)$ are also conformally equivalent, by the same isomorphism. Therefore if $e^{i\varphi_1},e^{i\varphi_2}$ are in the same subset $S_+$ or $S_-$, then $e^{i(\varphi_1+\varphi)},e^{i(\varphi_2+\varphi)}$ are also in the same subset for every $\varphi$. In particular, $1$ and $e^{i\varphi}$ are in the same subset if and only if $e^{-i\varphi}$ and $1$ are in the same subset. It follows that $e^{i\varphi},e^{-i\varphi}$ are in the same subset. But then $e^{2i\varphi},1$ are in the same subset, for every $\varphi$. Since $1 \in S_+$, we obtain $S_- = \emptyset$. This completes the proof.
\end{proof}

\subsection{Existence of the Killing vector field} \label{subs:Killing}

We now proceed to the second step, showing that domains with solutions of \eqref{Wangs_equation} corresponding to self-associated cones carry a Killing vector field with the properties claimed at the beginning of Section \ref{sec:Killing}. Note that in the complex isothermal coordinate $z$ on the domain $M$ a vector field on $M$ is represented by a complex-valued function $\psi$. 
A Killing vector field by definition generates a group of orientation-preserving isometries of $M$, which are in particular holomorphic functions on $M$. Then Killing vector fields are represented by holomorphic functions $\psi$ in any isothermal coordinate, since they are limits of sequences of holomorphic functions with uniform convergence on compact subsets.

{\lemma \label{lem:Killing} Let $K \subset \mathbb R^3$ be a self-associated cone, and let $(u,U)$ be a corresponding solution of Wang's equation \eqref{Wangs_equation} on a domain $M \subset \mathbb C$. Then there exists a holomorphic function $\psi$ on $M$ satisfying the relations
\begin{align} \label{KillingConditionU}
iU(z) + U'(z)\psi(z) + 3U(z)\psi'(z) &= 0, \\
\label{KillingConditionu}
Re(u'(z)\psi(z) + \psi'(z)) &= 0,
\end{align}
where the prime denotes the derivative with respect to the complex coordinate $z$. }

\begin{proof}
Assume the conditions of the lemma.

If $U \equiv 0$, then $\psi(z) \equiv 0$ satisfies the requirements. Let us henceforth assume that $U \not\equiv 0$.

By Fact \ref{fact:basic_equivalence} and Lemma \ref{lem:SLorbit} all solutions $(u,e^{i\varphi}U)$ of \eqref{Wangs_equation} on $M$ are mutually conformally equivalent. Let $\Aut M$ be the finite-dimensional Lie group of biholomorphic automorphisms of $M$. For every $e^{i\varphi} \in S^1$, let $G_{\varphi} \subset \Aut M$ be the set of all conformal isomorphisms taking the solution $(u,e^{i\varphi}U)$ of \eqref{Wangs_equation} to $(u,U)$. Since $e^{i\varphi}U \not\equiv e^{i\varphi'}U$ and hence $G_{\varphi} \cap G_{\varphi'} = \emptyset$ whenever $e^{i\varphi} \not= e^{i\varphi'}$, we may define a surjective map $\alpha: G = \bigcup_{e^{i\varphi} \in S^1} G_{\varphi} \to S^1$ by sending $G_{\varphi}$ to $e^{i\varphi}$.

Since the action of the automorphism group on the solution $(u,U)$ is continuous in the group element, the subgroups $G,G_0$ are closed in $\Aut M$, and hence $G_0$ is closed in the subspace topology of $G$. It follows that $G_0,G$ are Lie subgroups of $\Aut M$, and $G_0$ is a Lie subgroup of $G$. Since $G_0$ is the subgroup leaving $(u,U)$ and hence $(u,e^{i\varphi}U)$ invariant for every $\varphi$, the map $\alpha$ is a group homomorphism with kernel $G_0$.

If a sequence $b_k: w \mapsto z_k$ of automorphisms $b_k \in \Aut M$ tends to the identity map, then for every fixed $w \in M$ we have $b_k(w) \to w$ and $b'_k(w) \to 1$. Therefore $\alpha$ is continuous by \eqref{Wang_transform}. Thus $\alpha$ is a Lie group homomorphism and defines a group isomorphism between $S^1$ and the quotient $G/G_0$.

Let $\mathfrak g,\mathfrak g_0$ be the Lie algebras of the Lie groups $G,G_0$, respectively. Then $\alpha$ defines a Lie algebra homomorphism $\mathfrak a: \mathfrak g \to \mathbb R$ of $\mathfrak g$ into the Lie algebra $\mathbb R$ of $S^1$, such that $\mathfrak g_0 \subset \ker\mathfrak a$. For the sake of contradiction, let us assume that $\ker\mathfrak a = \mathfrak g$. Then the connection component of the identity in $G$ is mapped to $1 \in S^1$ and is hence a subset of $G_0$. It follows that $G_0$ is open in $G$, and the quotient topology of $G/G_0$ is discrete. However, $G/G_0$ is an uncountable set, because $S^1$ is uncountable. Hence $G$ is not second-countable, contradicting the property that it is a finite-dimensional Lie group. Thus $\ker\mathfrak a$ must be a strict Lie subalgebra of $\mathfrak g$, and $\mathfrak a$ is surjective.

Let $\mathfrak b \in \mathfrak g$ be such that $\mathfrak a(\mathfrak b) = -1$. Let $\{g_t\}_{t \in \mathbb R} \subset G$, $g_t: w \mapsto z_t$, be the one-parametric subgroup generated by the Lie algebra element $\mathfrak b$. Then $\alpha(g_t) = e^{-it}$, and from \eqref{Wang_transform} we get
\begin{equation} \label{KillingIntermediate}
e^{-it}U(w) = U(z_t)\left(\frac{dz_t}{dw}\right)^3,\qquad u(w) = u(z_t) + 2\log\left|\frac{dz_t}{dw}\right|
\end{equation}
for every $t \in \mathbb R$.

Let $\psi(w) = \left. \frac{dz_t(w)}{dt} \right|_{t = 0}$ be the velocity field of the flow generated by $\mathfrak b$, represented as a complex-valued function on $M$. Since the $g_t$ are isometries, the vector field $\psi$ is a Killing vector field and hence a holomorphic function. By virtue of the relations
\[ z_0(w) \equiv w,\qquad \left. \frac{dz_t}{dw} \right|_{t = 0} \equiv 1,\qquad \left[\frac{d}{dt}\frac{dz_t}{dw}\right]_{t = 0} = \psi'(w)
\]
differentiation of \eqref{KillingIntermediate} with respect to $t$ at $t = 0$ yields
\[ -iU(w) = U'(w)\psi(w) + 3U(w)\psi'(w),\qquad 0 = 2Re\,(u'(w)\psi(w)) + 2Re\,\psi'(w).
\]
The claim of the lemma readily follows.
\end{proof}

The existence of a Killing vector field $\psi$ satisfying \eqref{KillingConditionU},\eqref{KillingConditionu} is not only necessary, but also sufficient for the cone $K$ to be self-associated. We have the following result.

{\lemma \label{lem:reverse_Killing} Let $K \subset \mathbb R^3$ be a cone, and let $(u,U)$ be a solution of Wang's equation \eqref{Wangs_equation} corresponding to $K$, defined on a domain $M \subset \mathbb C$. Suppose there exists a holomorphic function $\psi$ on $M$ satisfying relations \eqref{KillingConditionU},\eqref{KillingConditionu}. Then $\psi$ generates a 1-parametric subgroup of conformal automorphisms of $M$, and the cone $K$ is self-associated. }

\begin{proof}
Let $\{g_t\}_{t \in \mathbb R}$ be the one-parametric family of conformal maps generated by the vector field $\psi$. For every point $w \in M$, $z_t = g_t(w)$ is defined at least in some neighbourhood of $t = 0$. Integrating \eqref{KillingConditionU},\eqref{KillingConditionu} along the trajectories of the flow we obtain \eqref{KillingIntermediate}. In particular, the maps $g_t$ are isometries, and the speed along the curve $(g_t(w))_t$ is constant for every fixed $w \in M$. By virtue of the completeness of the metric $h$ on $M$, $g_t$ is hence defined for all $t \in \mathbb R$ and is actually a conformal automorphism of $M$. 

Relation \eqref{KillingIntermediate} shows that $g_t$ takes the solution $(u,e^{-it}U)$ of \eqref{Wangs_equation} to the solution $(u,U)$ for every $t \in \mathbb R$. By Fact \ref{fact:basic_equivalence} all these solutions correspond to the same $SL(3,\mathbb R)$-orbit of cones, and $K$ is self-associated.
\end{proof}

Finally we show that it is actually sufficient to demand that $\psi$ generates a 1-parametric group of automorphisms of $M$ and satisfies condition \eqref{KillingConditionU} only.

{\lemma \label{lem:onlyU} Let $(u,U)$ be a complete solution of Wang's equation \eqref{Wangs_equation} on a simply connected domain $M \subset \mathbb C$. Let the holomorphic function $\psi$ generate a 1-parametric group of conformal automorphisms $g_t$ of $M$ and satisfy condition \eqref{KillingConditionU}. Then condition \eqref{KillingConditionu} holds. }

\begin{proof}
Fix $t \in \mathbb R$. Let $u_0$ be the complete solution of Wang's equation for $U$ and $u_t$ the complete solution for $e^{-it}U$ on $M$. By Fact \ref{fact:uniqueness} both exist and are unique. As in the proof of Lemma \ref{lem:reverse_Killing}, condition \eqref{KillingConditionU} implies that the automorphism $g_t$ transforms the holomorphic cubic differential $e^{-it}U$ to $U$, and hence also $u_t$ to $u_0$.

However, equation \eqref{Wangs_equation} is left invariant if $U$ is multiplied by a constant phase factor. Therefore the solutions $u_0,u_t$ actually coincide. It follows that $g_t$ preserves the affine metric $e^{u_0}|dz|^2$ for all $t$, and $\psi$ must be a Killing vector field. Condition \eqref{KillingConditionu} now follows as in the proof of Lemma \ref{lem:Killing}.
\end{proof}

\medskip

In this section we expressed the property of a cone $K \subset \mathbb R^3$ to be self-associated by an analytic condition on the corresponding cubic differential $U$. Namely, there must exist an infinitesimal generator $\psi$ of a 1-parametric subgroup of automorphisms of the domain $M$ of $U$ such that condition \eqref{KillingConditionU} holds. In the next section we shall classify the holomorphic functions $U$ together with their domains $M$ which satisfy condition \eqref{KillingConditionU}.

\section{Classification of cubic holomorphic differentials} \label{sec:cubic_class}

In this section we classify, up to biholomorphic isomorphisms, all pairs $(\psi,U)$ of holomorphic functions on $M$ satisfying condition \eqref{KillingConditionU}, where $\psi$ is representing a vector field generating a 1-parametric subgroup of biholomorphic automorphisms of $M$. At first we shall single out the following special case, however.

\medskip

\emph{Case 0:} $U \equiv 0$. This condition implies that the cubic form $C$ of the affine sphere vanishes identically. The Theorem of Pick and Berwald \cite[Theorem 4.5, p.~53]{NomizuSasaki} states that this happens if and only if the affine sphere is a quadric. A quadric can be asymptotic to the boundary of a cone in $\mathbb R^3$ only if this cone is ellipsoidal, i.e., linearly isomorphic to the Lorentz cone
\[ L_3 = \left\{ x = (x_0,x_1,x_2)^T \,|\, x_0 \geq \sqrt{x_1^2 + x_2^2} \right\}.
\]
The affine sphere is then one sheet of a two-sheeted hyperboloid and is isometric to a hyperbolic space form. The domain of the solution can be transformed to the unit disc $\mathbb D$ by a conformal isomorphism.


\medskip

In the rest of this section we assume that $U \not\equiv 0$. Then we also have $\psi \not\equiv 0$. The simply connected domains $M \subset \mathbb C$ admitting a 1-parameter group of biholomorphic automorphisms have been classified up to conformal isomorphisms and can be reduced to the following five cases \cite[Theorem 3.2, (2)]{DorfmeisterMa16a}, the domains in the other items of this theorem being either compact or not simply connected:
\begin{itemize}
\item[(a)] $M = \mathbb C$, $\psi(z) \equiv \mu$;
\item[(b)] $M = \mathbb C$, $\psi(z) = i\mu z$;
\item[(c)] $M = \mathbb D$, $\psi(z) = i\mu z$;
\item[(d)] $M = \{ z \,|\, Im\,z > 0 \}$, $\psi(z) \equiv \mu$;
\item[(e)] $M = \{ z \,|\, 0 < Im\,z < \pi \}$, $\psi(z) \equiv \mu$.
\end{itemize}
Here $\mu$ is an arbitrary non-zero real constant determining the parametrization of the automorphism subgroup. The cases (a),(d),(e) correspond to real translations, the cases (b),(c) to rotations about the origin. For each case we first find all solutions $U$ of \eqref{KillingConditionU}. In a second step we transform the domain by a conformal isomorphism to obtain a standardized form for $U$.

\medskip

\emph{Case T:} If the automorphism group consists of translations, the solutions of \eqref{KillingConditionU} are given by $U(z) = \gamma e^{-iz/\mu}$, where $\gamma$ is an arbitrary non-zero complex number. After application of the inverse of the conformal isomorphism $w \mapsto z = \alpha w + \beta$, where $\alpha = i\mu$, $\beta = -\mu(\frac{\pi}{2} + i\log(\gamma\mu^3))$, we obtain by \eqref{Wang_transform} the solution $(\psi,U) = (-i,e^w)$.

The applied map transforms the domain $M$ to the complex plane $\mathbb C$ in case (a), to the right half-plane $\{ w \,|\, Re\,w > \log|\gamma\mu^3| \}$ in case (b) when $\mu > 0$, to the left half-plane $\{ w \,|\, Re\,w < \log|\gamma\mu^3| \}$ in case (b) when $\mu < 0$, to the vertical strip $\{ w \,|\, \log|\gamma\mu^3| < Re\,w < \log|\gamma\mu^3| + \frac{\pi}{\mu} \}$ in case (c) if $\mu > 0$, and to the vertical strip $\{ w \,|\, \log|\gamma\mu^3| + \frac{\pi}{\mu} < Re\,w < \log|\gamma\mu^3| \}$ in case (c) when $\mu < 0$.

The solution of \eqref{KillingConditionU} can hence be brought to the form $(\psi,U) = (-i,e^w)$ on the domain $M_{(a,b)} = \{ w \,|\, a < Re\,w < b \}$, which is a vertical strip of finite or infinite width, parameterized by $-\infty \leq a < b \leq +\infty$ or equivalently by the set of non-empty open intervals $(a,b) \subset \mathbb R$.

\medskip

\emph{Case R:} If the automorphism group consists of rotations, the solutions are locally given by $U(z) = \gamma z^{-(3\mu+1)/\mu}$, where $\gamma$ is an arbitrary non-zero complex number. This solution is defined in a neighbourhood of $z = 0$ if and only if $-\frac{3\mu+1}{\mu}$ is a natural number, i.e., $\mu = -\frac{1}{k+3}$ for some $k \in \mathbb N$. After application of the inverse of the conformal isomorphism $w \mapsto z = \gamma^{\mu}w$ we obtain by \eqref{Wang_transform} the solution $(\psi,U) = (-\frac{iw}{k+3},w^k)$. The applied map transforms the domain $M$ to the complex plane $\mathbb C$ in case (b), and to the open disc $B_R = \{ w \,|\, |w| < R \}$ of radius $R = |\gamma|^{-\mu} \in (0,+\infty)$ in case (c).

The solution of \eqref{KillingConditionU} can therefore be brought to the form  $(\psi,U) = (-\frac{iw}{k+3},w^k)$, where $k \in \mathbb N$ is a discrete parameter, on the domain $B_R = \{ w \,|\, |w| < R \}$, which is an open disc of radius $R \in (0,+\infty]$.

\medskip

{\lemma \label{lem:Uclassification} Let $K \subset \mathbb R^3$ be a self-associated cone. Then the cubic differential $U$ from any solution of \eqref{Wangs_equation} corresponding to $K$ can be transformed by a conformal isomorphism to exactly one of the following cases:

\begin{itemize}
\item[0:] $U \equiv 0$ on $M = \mathbb D$;
\item[R:] $U = z^k$, $k \in \mathbb N$, on $M = B_R = \{ z \,|\, |z| < R \}$, $R \in (0,+\infty]$;
\item[T:] $U = e^z$ on $M = M_{(a,b)} = \{ z \,|\, a < Re\,z < b \}$, $-\infty \leq a < b \leq +\infty$.
\end{itemize}

Every listed case corresponds to an $SL(3,\mathbb R)$-orbit of self-associated cones, and the corresponding solution $u$ of \eqref{Wangs_equation} satisfies condition \eqref{KillingConditionu} with $\psi = 0$ in Case 0, $\psi = -\frac{iz}{k+3}$ in Case R, and $\psi = -i$ in Case T. }

\begin{proof}
Let $K \subset \mathbb R^3$ be a self-associated cone, and $(u,U)$ a complete solution of \eqref{Wangs_equation} corresponding to $K$. By Lemma \ref{lem:Killing} there exists a holomorphic function $\psi$ on $M$ satisfying \eqref{KillingConditionU},\eqref{KillingConditionu}. By Lemma \ref{lem:reverse_Killing} $\psi$ generates a 1-parameter group of conformal automorphisms of the domain $M$ of the solution. By the classification provided in this section the cubic differential can be transformed to exactly one of the cases in the lemma by a biholomorphic isomorphism of the domain. This proves the first assertion.

Let now $U: M \to \mathbb C$ be one of the holomorphic functions listed in the lemma. Then $U$ satisfies \eqref{KillingConditionU}, and the corresponding function $\psi$ generates a 1-parameter group of automorphisms of $M$. Let $u$ be the unique complete solution of \eqref{Wangs_equation}, which exists by Fact \ref{fact:uniqueness}. By Fact \ref{fact:basic_equivalence} this solution corresponds to an $SL(3,\mathbb R)$-orbit of cones. By Lemma \ref{lem:onlyU} the function $u$ satisfies \eqref{KillingConditionu}. But then $(u,U)$ corresponds to an $SL(3,\mathbb R)$-orbit of self-associated cones by Lemma \ref{lem:reverse_Killing}.
\end{proof}

Thus the zero case corresponds to a single $SL(3,\mathbb R)$-orbit, in the rotational case we obtain a countably infinite number of 1-parametric families of $SL(3,\mathbb R)$-orbits, while in the translational case we obtain a 2-parametric family of $SL(3,\mathbb R)$-orbits of self-associated cones.

\medskip

In this section we classified the cubic differentials $U$ together with their domain $M$ of definition corresponding to self-associated cones, up to conformal isomorphism. In the next section we find the corresponding solutions $u$ of \eqref{Wangs_equation}.

\section{Solving Wang's equation} \label{sec:solving_Wang}

In the previous section we have classified all solutions $(\psi,U)$ of equation \eqref{KillingConditionU}, which in the case $U \not\equiv 0$ led to the canonical forms $z^k$ and $e^z$ for the holomorphic function $U$ on families of domains $M \subset \mathbb C$. In this section we obtain the corresponding real-valued function $u$ by solving \eqref{Wangs_equation}. We use condition \eqref{KillingConditionu}, which holds by Lemma \ref{lem:Uclassification}, to reduce Wang's equation to an ODE. This ODE turns out to be equivalent to the Painlev\'e III equation, which is considered in more detail in the Appendix. The boundary conditions for the solution of the ODE are obtained from the completeness condition on the metric $h = e^u|dz|^2$ on $M$. Existence and uniqueness of the solution are guaranteed by Fact \ref{fact:uniqueness}.

%
%
%
%

\smallskip

We shall consider the Cases R and T from Lemma \ref{lem:Uclassification} separately.

\medskip

\emph{Case R:} Plugging $\psi = -\frac{iz}{k+3}$ into \eqref{KillingConditionu}, we obtain that $u$ is invariant under rotations of the domain $B_R$ about the origin, i.e., $u(z) = \chi(r)$ for some function $\chi: [0,R) \to \mathbb R$ which can be analytically extended to an even function on $(-R,R)$. Plugging $u(z) = \chi(r)$ into equation \eqref{Wangs_equation}, we obtain the second-order ODE
\begin{equation} \label{rotODE}
\frac{d^2\chi}{dr^2} = 2e^{\chi} - \frac{1}{r}\frac{d\chi}{dr} - 4r^{2k}e^{-2\chi}.
\end{equation}
Set $t_R = \sqrt{\frac{32R^{k+3}}{(k+3)^3}} \in (0,+\infty]$. Making the substitution $e^{\chi} = \sqrt{\frac{k+3}{2}}vr^{(k-1)/2}$, $t = \sqrt{\frac{32}{(k+3)^3}}r^{(k+3)/2}$ in \eqref{rotODE}, we obtain a solution $v(t)$ of the Painlev\'e equation \eqref{PainleveIII} which is positive on $(0,t_R)$. It hence equals one of the functions $v_{s,c}(t)$, $(s,c) \in [-1,3] \times \mathbb R$, from Proposition \ref{prop:Kitaev} in the Appendix. We now turn to the boundary conditions.

It is checked by straightforward calculation that the analyticity condition at $r = 0$ is equivalent to the condition that the corresponding Painlev\'e transcendent is of the form $v_{s,c}(t) = t^{-(k-1)/(k+3)}\cdot \kappa(t^{4/(k+3)})$, where $\kappa$ is an analytic function in some neighbourhood of the origin satisfying 
\[ \kappa(0) = e^{u(0)}2^{(3k-1)/(k+3)}(k+3)^{-2k/(k+3)} > 0.
\]
Comparing with the asymptotic expressions in Proposition \ref{prop:Kitaev}, we see that the parameter $s = s_k$ of the solution depends on $k$ only and is given by $\frac{s_k-1}{2} = \cos\frac{2\pi}{k+3}$, while the parameter $c$ is determined by the value
\[ \kappa(0) = \left\{ \begin{array}{rcl} e^c,&\quad& k = 0; \\ \lambda_k^28^{1-2\lambda_k}e^{3\lambda_kc/2}\frac{\Gamma(1-\frac{\lambda_k}{2})\Gamma(1-\lambda_k)}{(1+\frac{\lambda_k}{2})\Gamma(1+\lambda_k)},&& k > 0. \end{array} \right.
\]
Here $\lambda_k = \frac{2}{k+3}$.

Completeness of the metric $h = e^u|dz|^2$ on $B_R$ implies that the metric $e^{\chi}dr^2$ is complete on $(-R,R)$. In particular, in the case $R < \infty$ we must have $\lim_{r \to R}\chi(r) = +\infty$, which implies $\lim_{t \to t_R} v(t) = +\infty$. But then $v(t)$ must have a double pole with expansion \eqref{PIIIpole_expansion} at $t = t_R$.

{\lemma \label{lem:rot_unique} Let $R \in (0,\infty]$, $k \in \mathbb N$, $U(z) = z^k$ on $B_R$. Assume above definitions. Then the corresponding unique complete solution $u$ of \eqref{Wangs_equation} on $B_R$ is rotationally symmetric and determined by the unique positive solution $v_{s_k,c}(t)$ of the Painlev\'e III equation \eqref{PainleveIII} which is pole-free on $(0,t_R)$ and has a double pole at $t = t_R$ for $R$ finite, or is pole-free on the whole positive real axis for $R = +\infty$. If $R$ runs through the interval $(0,+\infty]$, then $c$ runs through $[0,+\infty)$ in the opposite direction. }

\begin{proof}
Existence and uniqueness of $u$ follows from Fact \ref{fact:uniqueness}, its rotational symmetry and expression through $v_{s_k,c}$ from the above considerations. 

For $R = +\infty$ the parameter $c$ must equal zero, because only for $c = 0$ the solutions $v_{s,c}$ are pole-free on a neighbourhood of $+\infty$ on the real axis. In particular, it follows that $v_{s_k,0}$ is pole-free and positive on $\mathbb R_+$, because it is the Painlev\'e transcendent which generates the complete solution $u$ on $\mathbb C$.

Let $R$ be finite. Note that $s_k \in (-1,3)$, and hence by Corollary \ref{cor:s_unique} there can exist at most one function $v_{s_k,c}$ which is positive on $(0,t_R)$ and has a double pole at $t = t_R$. Thus there exists a unique such function, namely that generating the solution $u$ on $B_R$.

Let now $c > 0$ be arbitrary. By Lemma \ref{lem:c_monotone} the solution $v_{s_k,c}$ is positive on $(0,t_c)$, where $t_c$ is the location of its left-most pole. But then by the preceding $v_{s_k,c}$ generates the solution $u$ of \eqref{Wangs_equation} on $B_R$, where $R$ is chosen such that $t_R = t_c$. For $c < 0$ the solution $v_{s_k,c}$ cannot be positive up to its left-most pole by Lemma \ref{lem:c_monotone}, and hence does not generate a complete solution $u$.

Hence the parameters $c \in [0,+\infty)$ and $R \in (0,+\infty]$ are in bijective correspondence. By Lemma \ref{lem:c_monotone} $c$ is a decreasing function of $R$, which proves the last assertion.
\end{proof}

%

\medskip

\emph{Case T:} 
%
%
Plugging $\psi = -i$ into \eqref{KillingConditionu}, we obtain that the solution $u$ is invariant with respect to vertical translations, i.e., of the form $u(z) = \chi(x)$ for some analytic function $\chi: (a,b) \to \mathbb R$. Plugging $u(z) = \chi(x)$ into equation \eqref{Wangs_equation}, we obtain the second-order ODE
\begin{equation} \label{translODE}
\frac{d^2\chi}{dx^2} = 2e^{\chi} - 4e^{2x}e^{-2\chi}.
\end{equation}
Set $t_a = \sqrt{32}e^{a/2} \in [0,+\infty)$, $t_b = \sqrt{32}e^{b/2} \in (0,+\infty]$. Making the substitution $e^{\chi} = \frac18vt$, $t = \sqrt{32}e^{x/2}$, we again obtain a solution $v(t)$ of the Painlev\'e equation \eqref{PainleveIII}.


Completeness of the metric $h = e^u|dz|^2$ on $M$ is equivalent to completeness of the metric $e^{\chi}dx^2$ on $(a,b)$. 

We shall treat the cases $a,b$ finite or infinite separately.

Consider first the case of finite $a$ ($b$). Then we must have $\lim_{x \to a}\chi(x) = +\infty$ ($\lim_{x \to b}\chi(x) = +\infty$). This condition is equivalent to $\lim_{t \to t_a} v(t) = +\infty$ ($\lim_{t \to t_b} v(t) = +\infty$). Hence $v(t)$ must have a double pole with expansion \eqref{PIIIpole_expansion} at $t = t_a$ ($t = t_b$). 


Let us now consider the case $a = -\infty$. Then the solution $v$ is pole-free on $(0,t_b)$ and hence as in Case R of the form $v_{s,c}$ for some $(s,c) \in [-1,3] \times \mathbb R$.
At the left end of the interval $(-\infty,b)$ the completeness condition is equivalent to the divergence of the integral $\int_{-\infty}^x e^{\chi(s)/2}\,ds = \int_0^t \sqrt{\frac{v(s)}{2s}}\,ds$. Comparing with the asymptotic expressions in Proposition \ref{prop:Kitaev} we see that this integral diverges if and only if $\lambda = 0$, or equivalently $s = 3$. 



In the case $b = +\infty$ the solution $v(t)$ is pole-free on the interval $(t_a,+\infty)$, and as in Case R above the parameter $c$ of the solution $v_{s,c}$ must equal zero.

{\lemma \label{lem:trans_unique} Let $-\infty \leq a < b \leq +\infty$, and let $U(z) = e^z$ be defined on $M_{(a,b)} = (a,b) + i\mathbb R$. Assume above definitions. Then the corresponding unique complete solution $u$ of \eqref{Wangs_equation} on $M_{(a,b)}$ is invariant with respect to vertical translations and determined by the unique positive solution $v(t)$ of the Painlev\'e equation \eqref{PainleveIII} which is pole-free on $(t_a,t_b)$ and
\begin{itemize}
\item of the form $v_{3,c}$ with a double pole at $t = t_b$ for $a = -\infty$, $b$ finite, here if $b$ runs through $\mathbb R$, then $c$ runs through $\mathbb R_+$ in the opposite direction;
\item of the form $v_{s,0}$ with a double pole at $t = t_a$ for $a$ finite, $b = +\infty$, here if $a$ runs through $\mathbb R$, then $s$ runs through $(3,+\infty)$ in the same direction;
\item with double poles at $t = t_a$, $t = t_b$ for $a,b$ finite;
\item equal to $v_{3,0}$ for $-a = b = +\infty$.
\end{itemize} }

\begin{proof}
Existence of $u$ follows from Fact \ref{fact:uniqueness}, its translational symmetry and expression through solutions of \eqref{PainleveIII} with the specified properties from the above considerations. 

The uniqueness of these solutions $v(t)$ follows from the fact that their asymptotics at $t = t_a$, $t = t_b$ implies completeness of the solution $u$ on $M_{(a,b)}$ constructed from them. But this complete solution is unique by Fact \ref{fact:uniqueness}.

Let us prove the assertions on the dependence of the parameters $s,c$ on $a,b$. Note that the solution $v_{3,0}$ is positive on $\mathbb R_+$, because it generates the complete solution $u$ of \eqref{Wangs_equation} on $\mathbb C$. 

By Lemma \ref{lem:s0monotone}, for every $s > 3$ the solution $v_{s,0}$ is positive on $(t_s,+\infty)$, where $t_s$ is the location of its right-most pole. Hence it generates the complete solution $u$ of \eqref{Wangs_equation} on $M_{(a,+\infty)}$, where $a$ is chosen such that $t_a = t_s$. For $s < 3$ the solution $v_{s,0}$ cannot have a pole on $\mathbb R$ and be positive up to the right-most pole by Lemma \ref{lem:s0monotone}, and hence does not generate a complete solution $u$. Therefore $s \in (3,+\infty)$ and $a \in \mathbb R$ are in bijective correspondence. By Lemma \ref{lem:s0monotone} $s$ is an increasing function of $a$.

Likewise, by Lemma \ref{lem:c_monotone}, for every $c > 0$ the solution $v_{3,c}$ is positive on $(0,t_s)$, where $t_s$ is the location of its left-most pole. Hence it generates the complete solution $u$ of \eqref{Wangs_equation} on $M_{(-\infty,b)}$, where $b$ is chosen such that $t_b = t_s$. For $c < 0$ the solution $v_{3,c}$ cannot have a pole on $\mathbb R$ and be positive up to the left-most pole by Lemma \ref{lem:c_monotone}, and hence does not generate a complete solution $u$. Therefore $c \in (0,+\infty)$ and $b \in \mathbb R$ are in bijective correspondence. By Lemma \ref{lem:c_monotone} $c$ is a decreasing function of $b$.
\end{proof}

\medskip


In this section we have reduced Wang's equation \eqref{Wangs_equation} to the Painlev\'e III equation \eqref{PainleveIII} and singled out those solutions of \eqref{PainleveIII} which generate the sought complete solutions of \eqref{Wangs_equation}. We are now ready to obtain the corresponding self-associated cones in the next sections.

\section{Integration of the frame equations} \label{sec:integr_frame}

In order to describe the cone $K$ corresponding to a solution $(u,U)$ of Wang's equation on a domain $M \subset \mathbb C$ we have to integrate the frame equations \eqref{frame_equations}. Employing the technique of \emph{osculation maps} introduced in \cite{DumasWolf15} we obtain a description of the boundary $\partial K$.

Recall that the affine sphere $f: M \to \mathbb R^3$ which is asymptotic to $\partial K$ is given by the third column of the unimodular moving frame $F: M \to SL(3,\mathbb C)$. In order to capture its asymptotics as the argument $z$ tends to the boundary $\partial M$ we compare $F$ to an explicit model frame $V: M \to SL(3,\mathbb C)$ (in \cite{DumasWolf15} the role of $V$ is played by the frame $F_T$ of a Titeica surface), such that the ratio $G = FV^{-1}$ (in \cite{DumasWolf15} the role of $G$ is played by the osculation map $\hat F$) remains finite at the boundary of the domain $M$. We can then read off $\partial K$ from the limit $G_0$ of $G$. Different pieces of the boundary necessitate different model frames, and in this section we compute only these individual pieces.

In order to treat Case R in the same framework as Case T, we apply the coordinate transformation 
\begin{equation} \label{rot_to_trans}
z \mapsto w = -3\log(k+3)+(k+3)\log z
\end{equation}
which transforms the cubic differential $U = z^k$ into $\tilde U = e^w$, the solution $\chi$ of \eqref{rotODE} into the solution $\tilde\chi = \chi - \frac{2k}{k+3}\log(k+3) + \frac{2x}{k+3}$ of \eqref{translODE}, the independent variable $r$ of \eqref{rotODE} into the independent variable $x = -3\log(k+3)+(k+3)\log r$ of \eqref{translODE}, and the punctured domain $B_R \setminus \{0\}$ into the (semi-)infinite horizontal strip $(-\infty,b) + i(-(k+3)\pi,(k+3)\pi]$ with $b = -3\log(k+3)+(k+3)\log R \in (-\infty,+\infty]$. It is more convenient, however, to work with the universal cover of $B_R \setminus \{0\}$, which is mapped to the left half-plane $(-\infty,b) + i\mathbb R$. We then have to keep in mind that the points $(x,y)$ and $(x,y+2(k+3)\pi)$ in this domain correspond to the same point on the affine sphere. In particular, the moving frame satisfies the periodicity relation
\begin{equation} \label{rot_periodicity}
F(x,y+2(k+3)\pi) = F(x,y)
\end{equation}
for all $x < b$, $y \in \mathbb R$. Since we are interested only in the behaviour of $F$ near the boundary of the domain $B_R$, the loss of the central point $z = 0$ is not relevant for the analysis.

In both Cases R and T the domain $M$ of definition of the solution $(u,U)$ is now given by $(a,b) + i\mathbb R$, where $-\infty \leq a < b \leq +\infty$, the cubic differential by $U = e^z$, and the conformal factor by $u(x,y) = \chi(x)$, where $\chi$ is a solution of ODE \eqref{translODE} on $(a,b)$.

The frame equations \eqref{frame_equations} take the form
\begin{equation} \label{equationsFexp}
\begin{aligned}
F_x &= F\begin{pmatrix} -e^{-\chi}e^x\cos y & e^{-\chi}e^x\sin y & e^{\chi/2} \\ e^{-\chi}e^x\sin y & e^{-\chi}e^x\cos y & 0 \\ e^{\chi/2} & 0 & 0 \end{pmatrix} =: FA,\\
F_y &= F\begin{pmatrix} e^{-\chi}e^x\sin y & -\frac{\chi'}{2} + e^{-\chi}e^x\cos y & 0 \\ \frac{\chi'}{2} + e^{-\chi}e^x\cos y & -e^{-\chi}e^x\sin y & e^{\chi/2} \\ 0 & e^{\chi/2} & 0 \end{pmatrix} =: FB.
\end{aligned}
\end{equation}
We shall now study the asymptotic behaviour of the frame $F(x,y)$ as $x$ tends to the limiting values $a$ or $b$ while $y$ remains fixed. For each value $y \in \mathbb R$, the last column $f(x,y)$ of $F(x,y)$ tends to a boundary ray of the cone $K$. We shall see that as $y$ changes, the limiting rays either sweep a boundary piece of the cone, or cluster at discrete values in case the boundary piece consists of just one ray or is polyhedral. The latter behaviour has been observed in \cite{DumasWolf15}. In case the boundary piece is smooth we shall describe it by a vector-valued ODE. Different initial conditions for this ODE lead to solutions which are related by a linear map. Thus an ODE is a convenient way to describe an $SL(3,\mathbb R)$-orbit of conic boundaries rather than the boundary of a single cone.

We shall treat the cases of finite and infinite limiting values separately, which yields 4 cases in total to consider.

\medskip

\emph{Case $x \to b = +\infty$} turns out to be directly covered by the results developed in \cite{DumasWolf15} and corresponds to polyhedral boundary pieces.

In order to apply the apparatus developed in \cite[Section 5]{DumasWolf15} we have to perform a coordinate transformation. Consider a horizontal strip $S = (a,+\infty) + i(y_0,y_0+3\pi) \subset M$ of width $3\pi$. Applying the coordinate transformation $z \mapsto w = 2^{1/6} \cdot 3 \cdot e^{z/3}$, we map $U = e^z$ to $\tilde U = \frac{\sqrt{2}}{2}$, the strip $S$ to some half-plane in $\mathbb C$ minus possibly a compact disc centered on the origin, vertical segments $x + i\cdot(y_0,y_0+3\pi)$ to semi-circles with radius $r = 2^{1/6} \cdot 3 \cdot e^{x/3}$, and $u(z) = \chi(x)$ to a radial function $\tilde u(w) = \tilde\chi(r) = \chi(x) - \frac13\log 2 - \frac{2x}{3}$.

Recall that $e^{\chi} = \frac{v_{s,0} \cdot t}{8}$, where $t = \sqrt{32}e^{x/2}$, and $v_{s,0}(t)$ is a solution of the Painlev\'e equation \eqref{PainleveIII}. From its asymptotics in Section \ref{sec:PainleveIII} we obtain that
\[ e^{\chi} \sim 2^{1/3}e^{2x/3} + 2^{-1/2}3^{-1/4}\pi^{-1/2}se^{x/2}e^{-2^{2/3}3^{3/2}e^{x/3}}
\]
as $x \to +\infty$ and hence
\[ e^{\tilde\chi} \sim 1 + 2^{-3/4}3^{1/4}\pi^{-1/2}sr^{-1/2}e^{-\sqrt{6}r}
\]
as $r \to \infty$. This is precisely the asymptotics obtained in Section 5 and needed for the results in Section 6 of \cite{DumasWolf15} to be applicable.

After the reverse transformation back to the coordinates $x,y$ on $M$ these results translate into the following limiting behaviour of the affine sphere $f$.

{\lemma \label{lem:binfty} Let $a \in [-\infty,+\infty)$, $M = (a,+\infty)+i\mathbb R$, and let the cubic differential $U = e^z$, $z = x + iy$, be defined on $M$. Let $(u,U)$, $u(z) = \chi(x)$, be the solution of Wang's equation \eqref{Wangs_equation} on $M$ which in Case R is obtained by transformation \eqref{rot_to_trans} from the (punctured) complete solution $(u,z^k)$ on $B_{+\infty} = \mathbb C$ (in which case $a = -\infty$), and in Case T is the complete solution on $M$. Let $f: M \to \mathbb R^3$ be a corresponding affine sphere, obtained by integration of the frame equations \eqref{equationsFexp}.

Then for $\frac{y}{2\pi}$ not an integer, the limit as $x \to +\infty$ of the ray through $f(x,y)$ depends only on the integer part of $\frac{y}{2\pi}$. The limits of two neighboring strips of width $2\pi$ are distinct. The ray through $f(x,y)$ with $\frac{y}{2\pi}$ an integer tends for $x \to +\infty$ to a proper convex combination of the limiting rays of the neighboring strips. The other convex combinations of these limiting rays can be obtained as the limit of the rays through $f(x,y(x))$ along a curved path $y = y(x)$ on $M$ as $x \to +\infty$. \qed }

{\corollary \label{cor:binfty} Assume the notations of Lemma \ref{lem:binfty}. The cone $K$ to whose boundary $f$ is asymptotic has a polyhedral boundary piece. The image under $f$ of the horizontal strip $(a,+\infty) + 2\pi i\cdot[n,n+1]$, $n \in \mathbb Z$, is a surface which is asymptotic to a conic wedge $W_n$, consisting of an extreme ray of $K$ and parts of the neighbouring two faces of $K$.

In Case R the union of the wedges $W_n$, $n \in \mathbb Z$, is the whole boundary $\partial K$, $K$ is polyhedral and has $k+3$ extreme rays. In Case T the boundary piece consisting of the wedges $W_n$ for all $n \in \mathbb Z$ is non-self-intersecting and has infinitely many extreme rays. }

\begin{proof}
The first part of the Corollary is a consequence of Lemma \ref{lem:binfty}. Let us prove the second part.

By virtue of \eqref{rot_periodicity}, in Case R we have $f(x,y) = f(x,y+2\pi(k+3)$ for all $(x,y) \in M$. Hence $W_n = W_{n+k+3}$ for all $n \in \mathbb Z$, and there are only $k+3$ distinct wedges. Moreover, the first and the last wedge in a sequence of $k+3$ consecutive wedges intersect in a boundary ray, and the boundary piece formed of these $k+3$ wedges closes in on itself. Hence the cone has $k+3$ extreme rays in total, and the boundary piece is already the whole boundary $\partial K$. That $K$ is polyhedral with $k+3$ extreme rays also follows from \cite[Theorem 6.3]{DumasWolf15}.

In Case T the rays through $f(x,y)$ for different points $(x,y) \in M$ are distinct, and the boundary piece $\bigcup_{n \in \mathbb Z} W_n$ is non-self-intersecting. It is composed of infinitely many wedges, and hence it contains infinitely many extreme rays.
\end{proof}

\medskip

\emph{Case $x \to a > -\infty$:} The solution $v(t)$ of the Painlev\'e equation \eqref{PainleveIII} corresponding to $\chi(x)$ has a double pole at $t_0 = \sqrt{32}e^{a/2}$ with expansion \eqref{PIIIpole_expansion}. Plugging this expansion into $e^{\chi} = \frac{vt}{8}$ and writing shorthand $\delta = x - a$, or equivalently, $t-t_0 = \sqrt{32}e^{a/2}(e^{\delta/2}-1)$, we obtain the expansions
\[ e^{\chi} = \delta^{-2} + \tilde\alpha + O(\delta^2), \quad e^{\chi/2} = \delta^{-1} + \frac{\tilde\alpha}{2}\delta + O(\delta^3),
\]
\[ e^{-\chi+\delta} = \delta^2 + \delta^3 + \left( \frac12-\tilde\alpha \right)\delta^4 + \left( \frac16-\tilde\alpha \right)\delta^5 + O(\delta^6),\quad \frac{\chi'}{2} = -\delta^{-1} + \tilde\alpha\delta + O(\delta^3),
\]
where $\tilde\alpha = \frac{\alpha e^{a/2}}{\sqrt{2}} - \frac{1}{48}$, and $\alpha$ is the parameter from \eqref{PIIIpole_expansion}.

Let us further consider the matrix variable $G = FV^{-1}$, where $V = \begin{pmatrix} 0 & 0 & \delta \\ 0 & 1 & 0 \\ -\delta^{-1} & 0 & \delta^{-1} \end{pmatrix}$. Since $V$ is unimodular, we have $G \in SL(3,\mathbb R)$. The frame equations \eqref{equationsFexp} become
\[ G_x = GVAV^{-1} - GV_xV^{-1},\qquad G_y = GVBV^{-1}.
\]
Inserting the above expansions, we obtain
\begin{eqnarray*}
G_x &=& G\begin{pmatrix}  \frac{\tilde\alpha}{2}\delta + O(\delta^3) &             0 &         - \delta - \frac{\tilde\alpha}{2}\delta^3 + O(\delta^5) \\
e^a(\delta+\delta^2)\sin y + O(\delta^3) &   e^a(\delta^2+\delta^3)\cos y + O(\delta^4) &     - e^a(\delta^3+\delta^4)\sin y + O(\delta^5) \\
e^a(1+\delta)\cos y + O(\delta^2) & - e^a(\delta+\delta^2)\sin y + O(\delta^3) & -\frac{\tilde\alpha}{2}\delta - e^a\delta^2\cos y + O(\delta^3) \end{pmatrix}, \\
G_y &=& G\begin{pmatrix}   0 &        1 + \frac{\tilde\alpha}{2}\delta^2 + O(\delta^4) &                              0 \\
\frac{3\tilde\alpha}{2} + e^a\delta\cos y + O(\delta^2) &  - e^a(\delta^2+\delta^3)\sin y + O(\delta^4) & 1 - \tilde\alpha\delta^2 + O(\delta^3) \\
    - e^a(1+\delta)\sin y + O(\delta^2) & \frac{3\tilde\alpha}{2} - e^a\delta\cos y + O(\delta^2) &    e^a(\delta^2+\delta^3)\sin y + O(\delta^4) \end{pmatrix}.
\end{eqnarray*}
The limit of the expression $G^{-1}\nabla G$ as $\delta \to 0$ (equivalently, $x \to a$) hence exists and the convergence is uniform in $y$. Therefore the function $G(x,y)$, which as a function of $x$ for every fixed $y$ is the solution of a linear ODE, extends analytically to some unimodular matrix-valued function $G_0(y) =: G(a,y)$. This function obeys the linear ODE
\begin{equation} \label{G0def_a_finite} 
\frac{dG_0}{dy} = G_0 \cdot \begin{pmatrix}   0 &        1  &                              0 \\
 \frac{3\tilde\alpha}{2}  &  0 & 1  \\
    - e^a\sin y  & \frac{3\tilde\alpha}{2}  &    0 \end{pmatrix}
\end{equation}
with the $2\pi$-periodic coefficient matrix being the limit of $G^{-1}G_y$ as $\delta \to 0$.

Denote the third column of $G_0$ by $f_0$. Equation \eqref{G0def_a_finite} allows to express the derivatives of $f_0$ as a function of $G_0$. In particular, we get
\begin{equation} \label{f_der_G}
(f_0'',f_0',f_0) = G_0 \cdot \begin{pmatrix} 1 & 0 & 0 \\ 0 & 1 & 0 \\ \frac{3\tilde\alpha}{2} & 0 & 1 \end{pmatrix}
\end{equation}
and hence 
\begin{equation} \label{det_f_01}
\det(f_0'',f_0',f_0) = 1.
\end{equation}
Hence the vector $f_0(y)$ never vanishes. Moreover, from \eqref{G0def_a_finite} it follows that it obeys the linear third-order ODE
\begin{equation} \label{transl_scalar_diffeq_left_f}
\frac{d^3f_0}{dy^3} - 3\tilde\alpha\cdot\frac{df_0}{dy} + e^a\sin y\,\cdot f_0 = 0
\end{equation}
with $2\pi$-periodic coefficients. We are now able to prove the following result.

{\lemma \label{lem:a_finite} Let $a \in \mathbb R$, $b \in (a,+\infty]$, $M = (a,b)+i\mathbb R$, and let the cubic differential $U = e^z$, $z = x + iy$, be defined on $M$. Let $(u,U)$, $u(z) = \chi(x)$, be the complete solution of Wang's equation \eqref{Wangs_equation} on $M$. Let $f: M \to \mathbb R^3$ be a corresponding affine sphere, obtained by integration of the frame equations \eqref{equationsFexp}, and let $K \subset \mathbb R^3$ be the cone to whose boundary $f$ is asymptotic.

Then for fixed $y \in \mathbb R$ the limit as $x \to a$ of the ray through $f(x,y)$ is a boundary ray $\rho_y$ of $K$. The union $\bigcup_{y \in \mathbb R} \rho_y$ of these boundary rays forms a non-self-intersecting boundary piece of $K$ which is analytic everywhere except the tip of the cone. It can be obtained as the conic hull of a vector-valued solution of ODE \eqref{transl_scalar_diffeq_left_f} on $\mathbb R$ satisfying \eqref{det_f_01}. }

\begin{proof}
Assume above notations. Fix $y \in \mathbb R$. The vector $f$ equals $G$ times the last column of $V$. Hence $\delta \cdot f(x,y) = G(x,y)\cdot ((x-a)^2,0,1)^T$, and this product tends to $f_0(y)$ for $x \to a$. Therefore the ray through $f(x,y)$ tends to the ray through $f_0(y)$ as $x \to a$. Since $f$ is asymptotic to $\partial K$, this limit $\rho_y$ is a boundary ray of $K$. By definition $f_0$ satisfies \eqref{transl_scalar_diffeq_left_f} and \eqref{det_f_01}.

Now by virtue of \eqref{det_f_01} the derivative $\frac{df_0}{dy}$ and the vector $f_0$ are linearly independent for every $y$. Hence locally for different $y$ the limit rays $\rho_y$ are distinct, the solution curve $f_0$ is analytic, and the boundary piece $\bigcup_{y \in \mathbb R} \rho_y$ generated by this curve is also analytic everywhere except at the origin. The rays through $f(x,y)$ for different points $(x,y) \in M$ are distinct, and hence the boundary piece is non-self-intersecting.
\end{proof}

\medskip

\emph{Case $x \to b < +\infty$:} Replacing $a$ by $b$ we obtain the same expansions and can perform the same constructions as in the previous case. However, since now $x \to b$ from the left, we have that $\delta = x - b$ is negative. Therefore the product $\delta \cdot f(x,y)$ will tend to a point on $-\partial K$ rather than on $\partial K$. In order to correct this sign change we shall multiply $f_0$ by $-1$, i.e., we define $f_0(y)$ as minus the third column of the matrix $G_0(y)$. Instead of \eqref{det_f_01} we then get
\begin{equation} \label{det_f_0m1}
\det(f_0'',f_0',f_0) = -1.
\end{equation}
However, the vector-valued function $f_0(y)$ is still a solution of the linear third-order ODE
\begin{equation} \label{transl_scalar_diffeq_right_f}
\frac{d^3f_0}{dy^3} - 3\tilde\alpha\cdot\frac{df_0}{dy} + e^b\sin y\,\cdot f_0 = 0
\end{equation}
with $\tilde\alpha = \frac{\alpha e^{b/2}}{\sqrt{2}} - \frac{1}{48}$, $\alpha$ being the parameter from expansion \eqref{PIIIpole_expansion} at $t_0 = \sqrt{32}e^{b/2}$. We get the following result.

{\lemma \label{lem:b_finite} Let $b \in \mathbb R$, $a \in [-\infty,b)$, $M = (a,b)+i\mathbb R$, and let the cubic differential $U = e^z$, $z = x + iy$, be defined on $M$. Let $(u,U)$, $u(z) = \chi(x)$, be the solution of Wang's equation \eqref{Wangs_equation} on $M$ which in Case R is obtained by transformation \eqref{rot_to_trans} from the (punctured) complete solution $(u,z^k)$ on $B_R$, where $b = -3\log(k+3)+(k+3)\log R$ (in which case $a = -\infty$), and in Case T is the complete solution on $M$. Let $f: M \to \mathbb R^3$ be a corresponding affine sphere, obtained by integration of the frame equations \eqref{equationsFexp}, and let $K \subset \mathbb R^3$ be the cone to whose boundary $f$ is asymptotic.

Then for fixed $y \in \mathbb R$ the limit as $x \to b$ of the ray through $f(x,y)$ is a boundary ray $\rho_y$ of $K$. The union $\bigcup_{y \in \mathbb R} \rho_y$ of these boundary rays forms a boundary piece of $K$ which is analytic everywhere except the tip of the cone. It can be obtained as the conic hull of a vector-valued solution of ODE \eqref{transl_scalar_diffeq_right_f} on $\mathbb R$ satisfying \eqref{det_f_0m1}.

In Case R the union $\bigcup_{y \in \mathbb R} \rho_y$ makes up the whole boundary $\partial K$, and the solution $f_0(y)$ is $2\pi(k+3)$-periodic. In Case T all rays $\rho_y$ are different, and the boundary piece is non-self-intersecting. }

\begin{proof}
Repeating the arguments in the proof of Lemma \ref{lem:a_finite}, we obtain that for fixed $y \in \mathbb R$ the ray through $f(x,y)$ tends to the ray through $f_0(y)$ as $x \to b$, and this limit ray is a boundary ray of $K$, denote it by $\rho_y$. Here $f_0$ is a vector-valued solution of ODE \eqref{transl_scalar_diffeq_right_f} satisfying \eqref{det_f_0m1}. Linear independence of $f_0,f_0'$ implies that locally the rays $\rho_y$ are different for different $y$, and the solution curve $f_0$ as well as the boundary piece generated by it are analytic except at the origin.

By virtue of \eqref{rot_periodicity}, in Case R we have $f(x,y) = f(x,y+2\pi(k+3)$ for all $(x,y) \in M$. But then $f_0(y) = f_0(y+2\pi(k+3))$, and the solution $f_0$ is $2\pi(k+3)$-periodic. Hence the boundary piece generated by $f_0$ closes in on itself, and makes up the whole boundary $\partial K$.

In Case T the rays through $f(x,y)$ are distinct for different points $(x,y)$, and the boundary piece generated by $f_0$ cannot close in on itself.
\end{proof}

\medskip

\emph{Case $x \to a = -\infty$:} Recall that $e^{\chi} = \frac{v_{3,c} \cdot t}{8}$, where $t = \sqrt{32}e^{x/2}$, and $v_{3,c}(t)$ is a solution of the Painlev\'e equation \eqref{PainleveIII}. From Lemma \ref{lem:asymp_v3c} we obtain the expansions
\[ e^{\chi} = \frac{1}{(-x+\alpha)^2} + O(e^{2x}x^2),\quad e^{-\chi+x} = (-x+\alpha)^2e^x + O(e^{3x}x^6),
\]
\[ e^{\chi/2} = \frac{1}{-x+\alpha} + O(e^{2x}x^3),\quad \frac{\chi'}{2} = \frac{1}{-x+\alpha} + O(e^{2x}x^3)
\]
as $x \to -\infty$, where $\alpha = \log 2 - 3\gamma - \frac32c$, and $\gamma$ is the Euler-Mascheroni constant.

Set $\delta = \frac{1}{-x+\alpha}$ and define the matrix-valued function $G = FV^{-1}$ with $V = \begin{pmatrix} \frac{\delta}{2} & 0 & \frac{\delta}{2} \\ 0 & 1 & 0 \\ -\delta^{-1} & 0 & \delta^{-1} \end{pmatrix}$. Again $V$ is unimodular and hence $G \in SL(3,\mathbb R)$. By virtue of the above expansions and the relation $\frac{d\delta}{dx} = \delta^2$ the frame equations become
\begin{align*}
G_x &= GVAV^{-1} - GV_xV^{-1} = G\begin{pmatrix} O(\delta^{-2}e^x) & O(\delta^{-1}e^x) & O(e^x) \\ O(\delta^{-3}e^x) & O(\delta^{-2}e^x) & O(\delta^{-1}e^x) \\ O(\delta^{-4}e^x) & O(\delta^{-3}e^x) & O(\delta^{-2}e^x) \end{pmatrix}, \\
G_y &= GVBV^{-1} = G\begin{pmatrix} O(\delta^{-2}e^x) & O(\delta^{-1}e^x) & O(e^x) \\ 2 + O(\delta^{-3}e^x) & O(\delta^{-2}e^x) & O(\delta^{-1}e^x) \\ O(\delta^{-4}e^x) & 2 + O(\delta^{-3}e^x) & O(\delta^{-2}e^x) \end{pmatrix}.
\end{align*}
Here the constants in the $O$ terms are uniformly bounded over $y \in \mathbb R$. It follows that $G(x,y)$ converges to some function $G_0(y)$ as $x \to -\infty$, and this function obeys
\begin{equation} \label{transl_left_inf}
\frac{dG_0}{dy} = G_0 \cdot \begin{pmatrix} 0 & 0 & 0 \\ 2 & 0 & 0 \\ 0 & 2 & 0 \end{pmatrix} \quad \Rightarrow \quad G_0(y) = G_0(0) \cdot \begin{pmatrix} 1 & 0 & 0 \\ 2y & 1 & 0 \\ 2y^2 & 2y & 1 \end{pmatrix}.
\end{equation}

{\lemma \label{lem:ainfty} Let $b \in (-\infty,+\infty]$, $M = (-\infty,b)+i\mathbb R$, and let the cubic differential $U = e^z$, $z = x + iy$, be defined on $M$. Let $(u,U)$, $u(z) = \chi(x)$, be the complete solution of Wang's equation \eqref{Wangs_equation} on $M$. Let $f: M \to \mathbb R^3$ be a corresponding affine sphere, obtained by integration of the frame equations \eqref{equationsFexp}, and let $K \subset \mathbb R^3$ be the cone to whose boundary $f$ is asymptotic.

Then for fixed $y \in \mathbb R$ the limit as $x \to -\infty$ of the ray through $f(x,y)$ is a boundary ray $\rho_0$ of $K$ which is independent of $y$. }

\begin{proof}
For fixed $y \in \mathbb R$ the product $\delta \cdot f(x,y) = G(x,y) \cdot (\frac{\delta^2}{2},0,1)^T$ tends to the last column of $G_0(y)$ for $x \to -\infty$, or equivalently $\delta \to 0$. By \eqref{transl_left_inf} the latter equals the last column of $G_0(0)$ and is hence a fixed vector, independent of $y$. It follows that the ray through $f(x,y)$ tends to the ray through $G_0(0)$, which is hence a boundary ray $\rho_0$ of the cone $K$.
\end{proof}

\medskip

In this section we described the boundary pieces of $K$ to which the affine sphere $f(x,y)$ is asymptotic as $x \to a$ or $x \to b$ for bounded $y$. In the cases $x \to \pm\infty$ we obtain a polyhedral boundary piece and a single ray, respectively, while for finite limits $x \to a$, $x \to b$ the boundary piece is the conic hull of a vector-valued solution of ODE \eqref{transl_scalar_diffeq_left_f},\eqref{transl_scalar_diffeq_right_f}, respectively. In the next section we assemble these pieces in order to construct the whole cones.

\section{Computing the self-associated cones} \label{sec:assemble}

In this section we describe the cones $K$ corresponding to the solutions $(u,U)$ analyzed in the previous sections. The most important tool to assemble the boundary pieces obtained in the previous section will be the automorphism group of $K$.

If the cubic differential $U$ is normalized to $e^z$, then the coefficient matrices $A,B$ in the frame equations \eqref{equationsFexp} are invariant with respect to the translation $(x,y) \mapsto (x,y+2\pi)$. Hence the gradient of the expression $F(x,y+2\pi)F^{-1}(x,y)$ vanishes, and there exists a unimodular matrix $T$ such that $F(x,y+2\pi) = TF(x,y)$ for all $x+iy \in M$. In particular, $f(x,y+2\pi) = Tf(x,y)$, and $T$ is a non-trivial linear automorphism of the cone $K$ to which the surface $f$ is asymptotic.


Let us introduce the diagonal matrix $D = \diag(1,-1,1)$. Then $A(x,-y) = DA(x,y)D$, $B(x,-y) = -DB(x,y)D$ for all $x+iy \in M$. Hence the gradient of the expression $F(x,-y)DF^{-1}(x,y)$ vanishes, and there exists a matrix $\Sigma$ with determinant equal to $-1$ such that $F(x,-y) = \Sigma F(x,y)D$ for all $x+iy \in M$. Again we have $f(x,-y) = \Sigma f(x,y)$, and $\Sigma$ is a linear automorphism of $K$. Clearly it satisfies the relation $\Sigma^2 = I$.

We have the relation $T\Sigma F(x,y) = TF(x,-y)D = F(x,-y+2\pi)D$ and hence
\[ (T\Sigma)^2F(x,y) = F(x,-(-y+2\pi)+2\pi)D^2 = F(x,y).
\]
It follows that $(T\Sigma)^2 = I$, and in particular, that $T$ and $T^{-1}$ are conjugated by $\Sigma$. Therefore $T$ must have spectrum $\{1,\lambda,\lambda^{-1}\}$ with $\lambda$ on the real line or on the unit circle.

{\definition \label{def:types} We shall call the self-associated cone $K$ of
\begin{itemize}
\item \emph{elliptic} type if its automorphism $T$ has spectrum $\{1,e^{i\varphi},e^{-i\varphi}\}$ with $\varphi$ not a multiple of $\pi$;
\item \emph{parabolic} type if $T$ has spectrum in $\{-1,+1\}$;
\item \emph{hyperbolic} type if $T$ has spectrum $\{1,\lambda,\lambda^{-1}\}$ with $|\lambda| > 1$.
\end{itemize}}

Besides the automorphisms $T,\Sigma$, the moving frame allows to define an invariant symmetric matrix. Set $J = \diag(1,-1,-1)$. Then $A(x,y-\pi)J + JA^T(-y) = 0$, $B(x,y-\pi)J - JB^T(-y) = 0$. Hence the gradient of the expression $\Omega = F(x,y-\pi)JF^T(x,-y)$ vanishes, and this unimodular matrix is constant over the domain $M$. If we replace $y$ by $-y+\pi$ we obtain the transpose of this matrix, which shows that $\Omega$ is symmetric.

\medskip

We shall now consider the cones corresponding to the different solutions $(u,U)$ studied in the previous sections.

{\lemma \label{lem:symR} In Case R the subgroup of automorphisms of $K$ generated by $T$ and $\Sigma$ is isomorphic to the dihedral group $D_{k+3}$. The cones $K$ are of elliptic type. }

\begin{proof}
It is more convenient to pass back to the domain $M = B_R$, with the cubic differential given by $U = z^k$. Let $F(z)$ be the moving frame of the corresponding affine sphere. Then $F(e^{2\pi i/(k+3)}z) = TF(z)$, $F(\bar z) = \Sigma F(z)D$ for all $z \in M$. In particular, $f(e^{2\pi i/(k+3)}z) = Tf(z)$, $f(\bar z) = \Sigma f(z)$.

Hence $f(0) \in \mathbb R^3$ is an eigenvector of both $T$ and $\Sigma$ with eigenvalue 1. The complementary invariant subspace is also shared by both $T$ and $\Sigma$ and is given by the tangent space to the surface $f(z)$ at $z = 0$. It can naturally be parameterized by the complex variable $z = x + iy \in \mathbb C$. By definition $T$ acts on this subspace by rotations $z \mapsto e^{2\pi i/(k+3)}z$ and $\Sigma$ by reflections $(x,y) \mapsto (x,-y)$. Hence the spectrum of $T$ is given by $\{1,e^{2\pi i/(k+3)},e^{-2\pi i/(k+3)}\}$ and the claims easily follow.
\end{proof}


\medskip

Let us now consider the cones corresponding to Case T.

{\lemma \label{lem:symT} In Case T the subgroup of automorphisms of $K$ generated by $T$ and $\Sigma$ is isomorphic to the infinite dihedral group $D_{\infty}$. The spectrum of $T$ is real. }

\begin{proof}
Suppose for the sake of contradiction that $T$ has a complex eigenvalue $e^{i\varphi}$. Then $T$ is diagonalisable, and the sequence $T^n$, $n \in \mathbb N$, accumulates to the identity matrix. But then the sequence of vectors $f(x,y+2\pi n)$ accumulates to $f(x,y)$, contradicting the fact that a complete hyperbolic affine sphere is an embedding. Hence the spectrum of $T$ is real.

The matrices $F(x,y+2\pi n)$, $n \in \mathbb Z$, are mutually distinct. Hence also $T^n$ are mutually distinct. Thus the group generated by $T,\Sigma$ is isomorphic to $D_{\infty}$.
\end{proof}

Therefore in Case T the corresponding cones are of either parabolic or hyperbolic type. We now establish that this depends on whether $a = -\infty$ or $a$ finite.

{\lemma \label{lem:a_inf_parabolic} Let $-\infty < b \leq \infty$ and let $K$ be the cone corresponding to the complete solution of Wang's equation on $(-\infty,b) + i\mathbb R$ with $U = e^z$. Then $K$ is of parabolic type, and $T$ has spectrum $\{1\}$ with a 1-dimensional eigenspace. The corresponding eigenvector generates a boundary ray $\rho_0$ of $K$. }

\begin{proof}
In the previous section the moving frame was represented as a product $F(x,y) = G(x,y)V(x)$, where $G,V$ are unimodular matrix-valued functions. The relation $F(x,y+2\pi) = TF(x,y)$ then yields the similar relation $G(x,y+2\pi) = TG(x,y)$. Passing to the limit $x \to -\infty$, we obtain the relation $G_0(y+2\pi) = TG_0(y)$ and by virtue of \eqref{transl_left_inf}
\[ T = G_0(0) \begin{pmatrix} 1 & 0 & 0 \\ 2(y+2\pi) & 1 & 0 \\ 2(y+2\pi)^2 & 2(y+2\pi) & 1 \end{pmatrix} \begin{pmatrix} 1 & 0 & 0 \\ 2y & 1 & 0 \\ 2y^2 & 2y & 1 \end{pmatrix}^{-1} G_0(0)^{-1} = G_0(0) \begin{pmatrix} 1 & 0 & 0 \\ 4\pi & 1 & 0 \\ 8\pi^2 & 4\pi & 1 \end{pmatrix} G_0(0)^{-1}.
\]
This shows that the only eigenvalue of $T$ is 1 and it has geometric multiplicity 1.

Moreover, the eigenvector of $T$ is given by the last column of $G_0(0)$. By construction of $G_0$ this vector lies on the boundary of $K$.
\end{proof}

{\corollary Assume the conditions of the previous lemma.

If $b < +\infty$, then the boundary of $K$ is given by the closure of the conic hull of a vector-valued solution of ODE \eqref{transl_scalar_diffeq_right_f} satisfying \eqref{det_f_0m1} and tending to $\rho_0$ both in the forward and the backward direction.

If $b = +\infty$, then the boundary of $K$ is given by the closure of an infinite chain of of 2-dimensional faces accumulating in both directions to $\rho_0$. }

\begin{proof}
Let $v \in \mathbb R^3$ be a non-zero vector. Then by Lemma \ref{lem:a_inf_parabolic} the sequence of rays generated by $T^nv$ tends to either $\rho_0$ or $-\rho_0$ as $n \to \pm\infty$. Clearly if $v$ lies on $\partial K$, this sequence also lies on $\partial K$ and hence tends to $\rho_0$.

Consider first the case $b < +\infty$. By Lemma \ref{lem:ainfty} $\partial K$ contains an analytic piece given by the conic hull of a vector-valued solution $f_0(y)$ of ODE \eqref{transl_scalar_diffeq_right_f} satisfying \eqref{det_f_0m1}. By construction we have $f_0(y+2\pi) = Tf_0(y)$ for all $y \in \mathbb R$. Setting $v = f_0(y_0)$ for $y_0$ arbitrary, we see that the ray through $f_0(y)$ tends to $\rho_0$ as $y \to \pm\infty$.

Let now $b = +\infty$. By Lemma \ref{lem:binfty} the boundary of $K$ contains an infinite chain of planar faces on which $T$ acts by construction by a shift. Setting $v$ equal to any non-zero point on this chain yields the desired conclusion.
\end{proof}

{\lemma Let $-\infty < a < b \leq +\infty$ and let $K$ be the cone corresponding to the complete solution of Wang's equation on $(a,b) + i\mathbb R$ with $U = e^z$. Then $K$ is of hyperbolic type, and $T$ has spectrum $\{1,\lambda,\lambda^{-1}\}$ with $\lambda > 1$. The boundary of $K$ is given by the closure of a union $\Upsilon^- \cup \Upsilon^+$ of two boundary pieces. Here $\Upsilon^-$ is an analytic boundary piece given by the conic hull of a vector-valued solution $f^-_0(y)$ of ODE \eqref{transl_scalar_diffeq_left_f} satisfying \eqref{det_f_01}. If $b < +\infty$, then $\Upsilon^+$ is also analytic and given by the conic hull of a vector-valued solution $f^+_0(y)$ of ODE \eqref{transl_scalar_diffeq_right_f} satisfying \eqref{det_f_0m1}. If $b = +\infty$, then $\Upsilon^+$ is a two-sided infinite chain of planar faces of $K$. }

\begin{proof}
Let us assume first that $b$ is finite. By Lemmas \ref{lem:a_finite}, \ref{lem:b_finite} the boundary $\partial K$ contains two analytic pieces $\Upsilon^{\pm}$ which are the conic hulls of vector-valued solutions $f^{\pm}_0(y)$ of ODEs \eqref{transl_scalar_diffeq_right_f},\eqref{transl_scalar_diffeq_left_f} satisfying \eqref{det_f_0m1},\eqref{det_f_01}, respectively. The automorphisms $T,\Sigma$ of $K$ act on these solutions by $Tf^{\pm}_0(y) = f^{\pm}_0(y+2\pi)$, $\Sigma f^{\pm}_0(y) = f^{\pm}_0(-y)$.

The rays through $f_0^+(y)$ tend to some boundary rays $\rho^+_{\pm}$ of $K$ as $y \to \pm\infty$, respectively. In the same way, the rays through $f_0^-(y)$ tend to boundary rays $\rho^-_{\pm}$. Clearly these limit rays must be generated by eigenvectors of $T$ corresponding to positive eigenvalues. Moreover, $\Sigma[\rho^{\pm}_{\pm}] = \rho^{\pm}_{\mp}$.

If all these four rays are mutually distinct, then any three of them cannot be coplanar, otherwise at least one of the analytic boundary pieces has to lie in a flat face of $K$, contradicting \eqref{det_f_01} or \eqref{det_f_0m1}. It follows that $T$ has an eigenspace of dimension 3 and must equal the identity, leading to a contradiction.

Hence we have that $\rho^-_+ = \rho^+_+ =: \rho_+$, $\rho^-_- = \rho^+_- =: \rho_-$, because these two relations can hold only simultaneously by virtue of the symmetry $\Sigma$. Since the two pieces $\Upsilon^{\pm}$ link together at both ends, they together with the limit rays $\rho_{\pm}$ make up the whole boundary $\partial K$.

Moreover, the product of the eigenvalues of $T$ corresponding to $\rho_{\pm}$ equals 1, because $\Sigma$ conjugates $T$ to $T^{-1}$.

Now suppose for the sake of contradiction that the eigenvectors of $T$ generating $\rho_{\pm}$ have eigenvalue 1. Then the whole 2-dimensional subspace spanned by $\rho_{\pm}$ is left fixed by $T$. However, this subspace intersects the interior of the cone $K$, contradicting the fact that the affine sphere $f(x,y)$ does not contain fixed points of $T$. Hence the eigenvalues corresponding to $\rho_{\pm}$ equal $\lambda^{\pm1}$ for some $\lambda > 1$, proving our claim.

In the case $b = +\infty$ the argument is similar, with reference to Lemma \ref{lem:binfty} instead of Lemma \ref{lem:b_finite} and with the analytic boundary piece generated by the solution $f^+_0(y)$ replaced by an infinite chain of 2-dimensional faces.
\end{proof}

For each of the pieces $\Upsilon^{\pm}$ we have a unimodular initial condition to choose when integrating the moving frame. Only one of these conditions can be chosen freely, corresponding to the choice of the representative in the $SL(3,\mathbb R)$-orbit of the cone. In order to choose the second initial condition correctly without having to compute the solution $v(t)$ of the Painlev\'e equation we will make use of the automorphisms $T,\Sigma$, and the invariant symmetric matrix $\Omega$ which are common for both boundary pieces.

Consider first the case $b < +\infty$. For any of the solutions $f_0(y)$ of ODE \eqref{transl_scalar_diffeq_left_f} or \eqref{transl_scalar_diffeq_right_f} defining a boundary piece, set $\Phi = (f_0'',f_0',f_0)$. By \eqref{f_der_G} we obtain
\[ G_0 = \Phi \begin{pmatrix} 1 & 0 & 0 \\ 0 & 1 & 0 \\ -\frac{3\tilde\alpha}{2} & 0 & 1 \end{pmatrix}.
\]
Moreover, by $f_0(y+2\pi) = Tf_0(y)$, $f_0(-y) = \Sigma f_0(y)$ we get $\Phi(y+2\pi) = T\Phi(y)$, $\Phi(-y) = \Sigma\Phi(y)D$. Note further that
\[ \lim_{\delta\to0} V(\delta)JV^T(\delta) = \lim_{\delta\to0} \begin{pmatrix} 0 & 0 & \delta \\ 0 & 1 & 0 \\ -\delta^{-1} & 0 & \delta^{-1} \end{pmatrix} \begin{pmatrix} 1 & 0 & 0 \\ 0 & -1 & 0 \\ 0 & 0 & -1 \end{pmatrix} \begin{pmatrix} 0 & 0 & -\delta^{-1} \\ 0 & 1 & 0 \\ \delta & 0 & \delta^{-1} \end{pmatrix} = \begin{pmatrix} 0 & 0 & -1 \\ 0 & -1 & 0 \\ -1 & 0 & 0 \end{pmatrix}.
\]
Hence
\begin{align*}
\Omega &= G(x,y-\pi)V(\delta)JV^T(\delta)G^T(x,-y) = G_0(y-\pi)\begin{pmatrix} 0 & 0 & -1 \\ 0 & -1 & 0 \\ -1 & 0 & 0 \end{pmatrix}G^T_0(-y) \\ =& \Phi(y-\pi)\begin{pmatrix} 1 & 0 & 0 \\ 0 & 1 & 0 \\ -\frac{3\tilde\alpha}{2} & 0 & 1 \end{pmatrix}\begin{pmatrix} 0 & 0 & -1 \\ 0 & -1 & 0 \\ -1 & 0 & 0 \end{pmatrix}\begin{pmatrix} 1 & 0 & -\frac{3\tilde\alpha}{2} \\ 0 & 1 & 0 \\ 0 & 0 & 1 \end{pmatrix}\Phi^T(-y) = \Phi(y-\pi)\begin{pmatrix} 0 & 0 & -1 \\ 0 & -1 & 0 \\ -1 & 0 & 3\tilde\alpha \end{pmatrix}\Phi^T(-y).
\end{align*}
Let now $\Phi_{\pm}(y)$ be the matrices corresponding to the solutions $f_0^{\pm}(y)$ defining the two pieces $\Upsilon^{\pm}$. Then we have
\begin{equation} \label{comp_init_cond}
\begin{aligned}
\Phi_+(2\pi)\Phi_+^{-1}(0) = \Phi_-(2\pi)\Phi_-^{-1}(0),\ &\ \Phi_+(0)D\Phi_+^{-1}(0) = \Phi_-(0)D\Phi_-^{-1}(0),\\ \Phi_+(-\pi)\begin{pmatrix} 0 & 0 & -1 \\ 0 & -1 & 0 \\ -1 & 0 & 3\tilde\alpha_b \end{pmatrix}\Phi_+^T(0) &= \Phi_-(-\pi)\begin{pmatrix} 0 & 0 & -1 \\ 0 & -1 & 0 \\ -1 & 0 & 3\tilde\alpha_a \end{pmatrix}\Phi_-^T(0),
\end{aligned}
\end{equation}
where $\tilde\alpha_a,\tilde\alpha_b$ are the constants in ODEs \eqref{transl_scalar_diffeq_left_f},\eqref{transl_scalar_diffeq_right_f}, respectively. Together with the conditions $\det\Phi_{\pm} = \pm1$ and the condition that the conic hulls of the solutions $f_0$ meet at the same rays $\rho_{\pm}$ as $y \to \pm\infty$, respectively, these relations define $\Phi_+$ uniquely for given $\Phi_-$ and vice versa.

Now we consider the case $b = +\infty$. The formulas involving $\Phi_-$ remain the same, but we have to find a suitable replacement for those involving $\Phi_+$. Recall that the ray through $f(x,(2n-1)\pi)$, $n \in \mathbb Z$, tends to an extreme ray $\rho_n$ of $K$, and that $T[\rho_n] = \rho_{n+1}$, $\Sigma[\rho_n] = \rho_{1-n}$ for all $n$. Choose a non-zero point $v_0 \in \rho_0$ and define $v_n = T^nv_0 \in \rho_n$ for all $n$. By the relation $(T\Sigma)^2 = I$ we then have $\Sigma v_n = v_{1-n}$. Application of the methods in \cite[Section 6]{DumasWolf15} leads to the formula $\Omega = |\det H|^{-2/3}H\begin{pmatrix} -1 & 0 & 0 \\ 0 & s & -1 \\ 0 & -1 & 0 \end{pmatrix}H^T$, where $H = (v_1,v_0,v_{-1})$ and $s$ is the parameter of the solution $v_{s,0}$ of the Painlev\'e equation \eqref{PainleveIII} determining the conformal factor $u$. It can be computed from $T$ by $s = \tr\,T$. If we normalize $v_0$ such that $\det H = -1$ and demand that $v_n \to \rho_{\pm}$ as $n \to \pm\infty$, then $v_0$ can be determined uniquely from the formula
\begin{equation} \label{comp_init_cond_inf}
H\begin{pmatrix} -1 & 0 & 0 \\ 0 & s & -1 \\ 0 & -1 & 0 \end{pmatrix}H^T = \Phi_-(-\pi)\begin{pmatrix} 0 & 0 & -1 \\ 0 & -1 & 0 \\ -1 & 0 & 3\tilde\alpha_a \end{pmatrix}\Phi_-^T(0).
\end{equation}

\medskip

In this section we described the boundaries of the self-associated cones corresponding to every solution $(u,U)$ of Wang's equation which has been obtained and analyzed in the previous sections. In the next section we summarize our findings in formal theorems.

\section{Results} \label{sec:results}

The self-associated cones can be grouped into three types according to Definition \ref{def:types}, which are described in the theorems below. Besides these, the only self-associated cones are the ellipsoidal cones, which form the $SL(3,\mathbb R)$-orbit of the Lorentz cone $L_3$ described in Case 0 of Section \ref{sec:cubic_class}.

{\theorem \label{thm:ell} The $SL(3,\mathbb R)$-orbits of self-associated cones of elliptic type are parameterized by a discrete parameter $k \in \mathbb N$ and a continuous parameter $R \in (0,+\infty]$. Every such cone possesses a subgroup of linear automorphisms which is isomorphic to the dihedral group $D_{k+3}$. This subgroup is generated by an automorphism $T$ with spectrum $\{1,e^{2\pi i/(k+3)},e^{-2\pi i/(k+3)}\}$ and a reflection $\Sigma$. The cone possesses an interior ray which is left fixed by every automorphism in the subgroup.

The complete hyperbolic affine 2-sphere which is asymptotic to a cone of elliptic type corresponding to the parameters $(k,R)$ possesses an isothermal parametrization with domain $B_R \subset \mathbb C$, the open disc of radius $R$, and cubic differential $U = z^k$. The corresponding affine metric $h = e^u|dz|^2$ is given by the conformal factor $e^{u(z)} = \sqrt{\frac{k+3}{2}}|z|^{(k-1)/2}v_{s_k,c}(t)$, where $t = \sqrt{\frac{32}{(k+3)^3}}|z|^{(k+3)/2}$, $s_k = 1 + 2\cos\frac{2\pi}{k+3}$, and $v_{s_k,c}$ is a solution of the Painlev\'e III equation \eqref{PainleveIII}. Here for every fixed $k$ the parameter $c$ is strictly monotonely decreasing in $R$ such that $c = 0$ for $R = +\infty$, and $c \to +\infty$ for $R \to 0$. For finite $R$ it is determined by the condition that $v_{s_k,c}$ is positive on $(0,t_0)$ and has a double pole at $t_0 = \sqrt{\frac{32}{(k+3)^3}}R^{(k+3)/2}$ with expansion \eqref{PIIIpole_expansion}.

For every $k \in \mathbb N$ the self-associated cones corresponding to the parameters $(k,+\infty)$ are linearly isomorphic to a cone over a regular $(k+3)$-gon.

The self-associated cones which correspond to the parameters $(k,R)$ for finite $R$ can be constructed as follows. Choose $W \in SL(3,\mathbb R)$ arbitrarily and solve ODE \eqref{transl_scalar_diffeq_right_f} with initial values $(f_0'',f_0',f_0)|_{y = 0} = -W$, where $e^b = \frac{R^{k+3}}{(k+3)^3}$ and $\tilde\alpha = \frac{\alpha t_0}{8} - \frac{1}{48}$. Then the solution $f_0: \mathbb R \to \mathbb R^3$ traces a closed analytic $2(k+3)\pi$-periodic curve in $\mathbb R^3$. The boundary of the cone is analytic at every non-zero point and can be obtained as the conic hull of this curve. The $SL(3,\mathbb R)$-orbit of self-associated cones corresponding to $(k,R)$ can be parameterized by the initial condition $W$. }

\begin{figure}
\centering
\includegraphics[width = 14.22cm,height = 9.22cm]{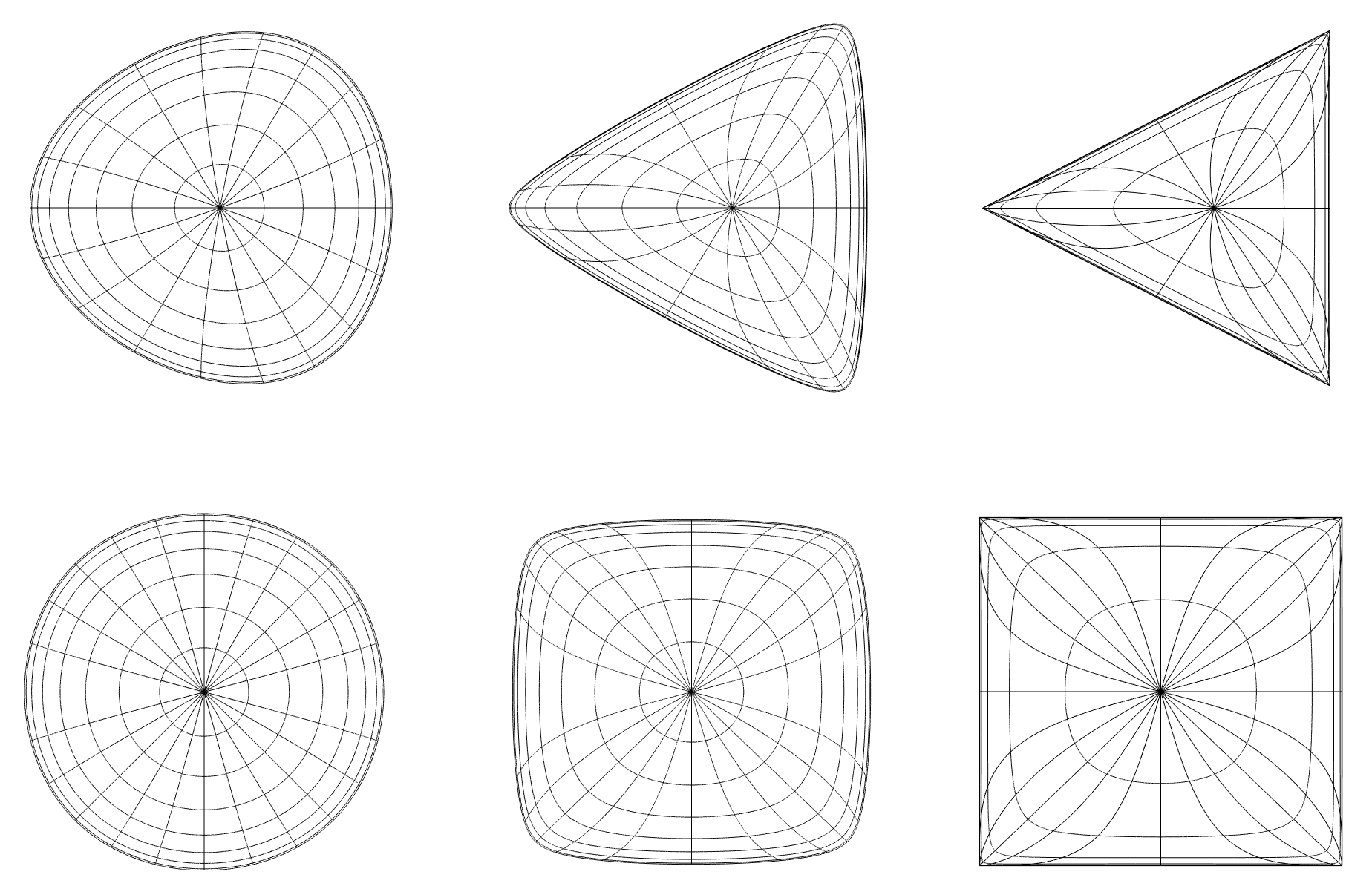}
\caption{Compact affine sections of self-associated cones of elliptic type for $k = 0,1$ and $R = 1,2,4$. The domain $M = B_R$ with uniform grid in polar coordinates is projected onto the interior of the section. The step size in the radial direction equals 1, in the angular direction $\frac{\pi}{3(k+3)}$.}
\label{fig_ell}
\end{figure}

In Fig.~\ref{fig_ell} we present compact affine sections of some self-associated cones of elliptic type. Since each interior ray intersects the affine sphere inscribed in the cone in exactly one point, we may project the isothermal coordinate $z \in B_R$ onto the interior of the affine section. A uniformly spaced polar coordinate grid on $B_R$ projects on the sections as shown in the figure.

{\theorem \label{thm:par} The $SL(3,\mathbb R)$-orbits of self-associated cones of parabolic type are parameterized by a continuous parameter $b \in (-\infty,+\infty]$. Every such cone $K$ possesses a subgroup of linear automorphisms which is isomorphic to the infinite dihedral group $D_{\infty}$. This subgroup is generated by an automorphism $T$ with spectrum $\{1\}$ and eigenvalue of geometric multiplicity 1, and a reflection $\Sigma$. There exists a boundary ray $\rho_0$ of $K$ which is left fixed by every automorphism in the subgroup.

The complete hyperbolic affine 2-sphere which is asymptotic to a cone of parabolic type corresponding to the parameter $b$ possesses an isothermal parametrization with domain $(-\infty,b) + i\mathbb R$ and cubic differential $U = e^z$. The corresponding affine metric $h = e^u|dz|^2$ is given by the conformal factor $e^{u(z)} = \frac{tv_{3,c}(t)}{8}$, where $t = \sqrt{32}e^{Re\,z/2}$, and $v_{3,c}$ is a solution of the Painlev\'e III equation \eqref{PainleveIII}. Here $c$ is a strictly monotonely decreasing function of $b$ such that $c = 0$ for $b = +\infty$, and $c \to +\infty$ for $b \to -\infty$. For finite $b$ it is determined by the condition that $v_{3,c}$ is positive on $(0,t_0)$ and has a double pole at $t_0 = \sqrt{32}e^{b/2}$ with expansion \eqref{PIIIpole_expansion}.

Any cone corresponding to the parameter value $b = +\infty$ is linearly isomorphic to the cone given by the closed convex conic hull of the infinite set of vectors $v_n = (n^2,n,1)^T$, $n \in \mathbb Z$.

The self-associated cones which correspond to a finite parameter $b$ can be constructed as follows. Choose $W \in SL(3,\mathbb R)$ arbitrarily and solve the vector-valued ODE \eqref{transl_scalar_diffeq_right_f} with initial values $(f_0'',f_0',f_0)|_{y = 0} = -W$, where $\tilde\alpha = \frac{\alpha e^{b/2}}{\sqrt{2}} - \frac{1}{48}$. Then the solution $f_0: \mathbb R \to \mathbb R^3$ traces an analytic curve in $\mathbb R^3$. The conic hull of this curve is analytic at every non-zero point and meets itself at the ray $\rho_0$. The boundary of the cone can be obtained as the union of the conic hull of $f_0$ with the ray $\rho_0$. The $SL(3,\mathbb R)$-orbit of self-associated cones corresponding to the parameter $b$ can be parameterized by the initial condition $W$. }

\begin{figure}
\centering
\includegraphics[width = 15.00cm,height = 14.50cm]{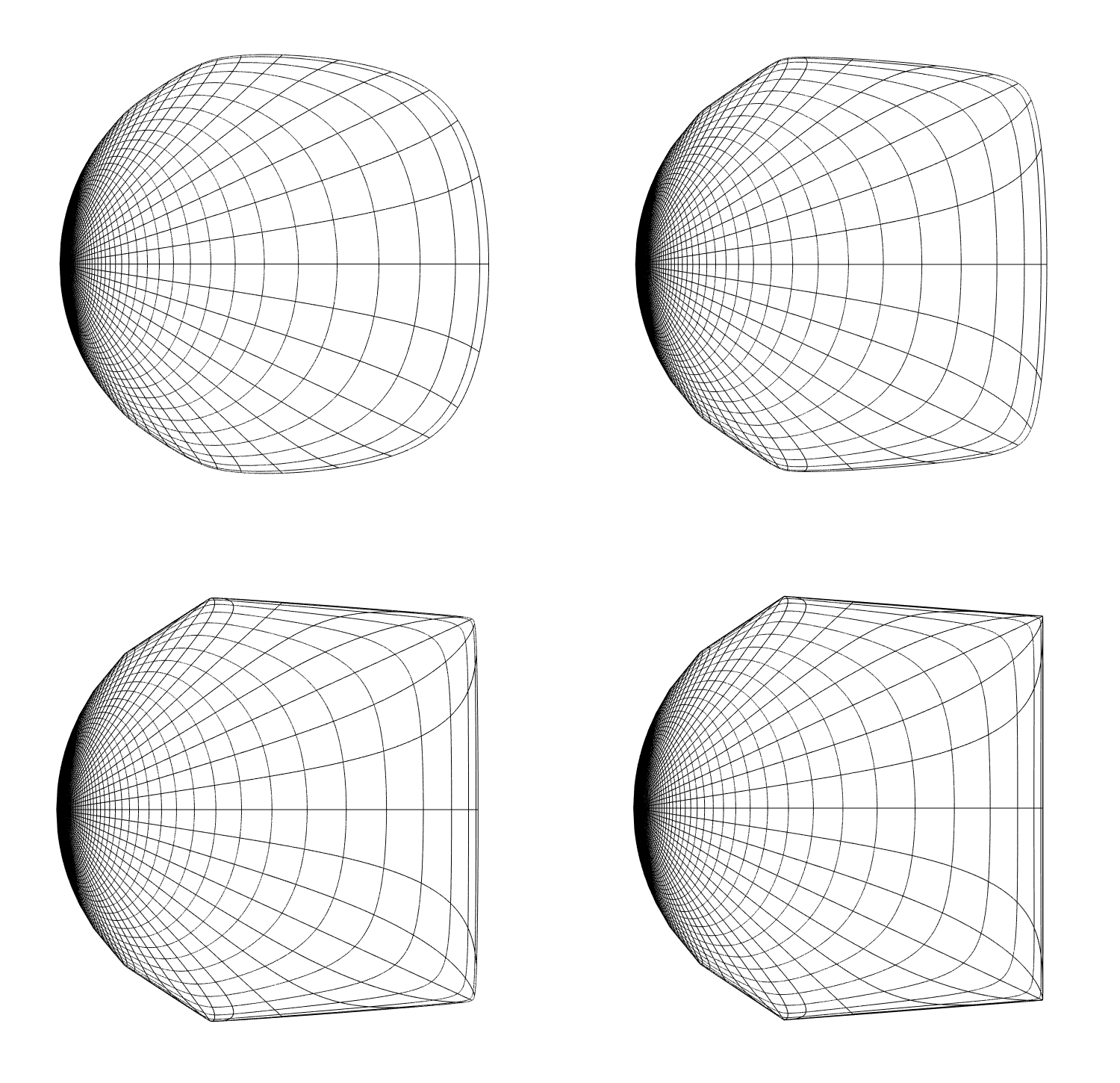}
\caption{Compact affine sections of self-associated cones of parabolic type for $b = -2,-1,0,1$. The domain $M = (-\infty,b) + i\mathbb R$ with uniform grid in cartesian coordinates is projected onto the interior of the section. The step size in the horizontal direction equals 1, in the vertical direction $\frac{\pi}{3}$. }
\label{fig_par}
\end{figure}

In Fig.~\ref{fig_par} we represent compact affine sections of self-associated cones of parabolic type for the parameter values $b = -2,-1,0,1$. As in Fig.~\ref{fig_ell}, we project a uniformly spaced grid in the domain $M$ on the section.

\medskip

{\theorem \label{thm:hyp} The $SL(3,\mathbb R)$-orbits of self-associated cones of hyperbolic type are parameterized by two continuous parameters $-\infty < a < b \leq +\infty$. Every such cone $K$ possesses a subgroup of linear automorphisms which is isomorphic to the infinite dihedral group $D_{\infty}$. This subgroup is generated by an automorphism $T$ with spectrum $\{1,\lambda,\lambda^{-1}\}$, $\lambda > 1$, and a reflection $\Sigma$. The cone possesses two boundary rays $\rho_{\pm}$ which are generated by eigenvectors of $T$ with eigenvalues $\lambda^{\pm1}$, respectively, and are mapped to each other by $\Sigma$.

The complete hyperbolic affine 2-sphere which is asymptotic to a cone of hyperbolic type corresponding to the parameters $(a,b)$ possesses an isothermal parametrization with domain $(a,b) + i\mathbb R$ and cubic differential $U = e^z$. The corresponding affine metric $h = e^u|dz|^2$ is given by the conformal factor $e^{u(z)} = \frac{tv(t)}{8}$, where $t = \sqrt{32}e^{Re\,z/2}$, and $v$ is a solution of the Painlev\'e III equation \eqref{PainleveIII}. This solution is characterized by the condition that it is positive on $(t_a,t_b)$ with $t_a = \sqrt{32}e^{a/2}$, $t_b = \sqrt{32}e^{b/2}$, has a double pole at $t_a$ with expansion \eqref{PIIIpole_expansion} featuring a constant $\alpha_a$, and for finite $b$ has a double pole at $t_b$ with expansion \eqref{PIIIpole_expansion} featuring a constant $\alpha_b$. If $b = +\infty$, then $v$ is given by the solution $v_{s,0}$, where $s$ is a strictly monotonely increasing function of $a$ such that $s \to 3$ for $a \to -\infty$, and $s \to +\infty$ for $a \to +\infty$.

The boundary of $K$ is the closure of the union of two pieces $\Upsilon^{\pm}$ which join at the rays $\rho_{\pm}$.

The piece $\Upsilon^-$ can be constructed as follows. Choose $W \in SL(3,\mathbb R)$ arbitrarily and solve the vector-valued ODE \eqref{transl_scalar_diffeq_left_f} with initial values $\Phi_-(0) = (f_0'',f_0',f_0)|_{y = 0} = W$, where $\tilde\alpha = \frac{\alpha_a e^{a/2}}{\sqrt{2}} - \frac{1}{48}$. Then the solution $f_0: \mathbb R \to \mathbb R^3$ traces an analytic curve in $\mathbb R^3$, whose conic hull is analytic at every non-zero point and equals $\Upsilon^-$. The curve $f_0$ tends to $\rho_{\pm}$ for $y \to \pm\infty$, and $T = (f_0'',f_0',f_0)|_{y = 2\pi} \cdot W^{-1}$.

The piece $\Upsilon^+$ can be constructed as follows.

\smallskip

Case $b < +\infty$: Solve the vector-valued ODE \eqref{transl_scalar_diffeq_right_f} with initial values $\Phi_+(0) = (f_0'',f_0',f_0)|_{y = 0} = -W'$, where $\tilde\alpha = \frac{\alpha_b e^{b/2}}{\sqrt{2}} - \frac{1}{48}$, and $W' \in SL(3,\mathbb R)$ is determined by relations \eqref{comp_init_cond} and the condition that $f_0(y)$ tends to $\rho_{\pm}$ as $y \to \pm\infty$, respectively. Then the solution $f_0: \mathbb R \to \mathbb R^3$ traces an analytic curve in $\mathbb R^3$, whose conic hull is analytic at every non-zero point and equals $\Upsilon^+$.

\smallskip

Case $b = +\infty$: Set $s = \tr\,T$ and let $v_0 \in \mathbb R^3$, $v_n = T^nv_0$, $n \in \mathbb Z$, be non-zero vectors such that \eqref{comp_init_cond_inf} holds with $H = (v_1,v_0,v_{-1})$, $\Sigma v_n = v_{1-n}$, and $v_n \to \rho_{\pm}$ for $n \to \pm\infty$, respectively. Then the right piece of $\partial K$ is given by a two-sided infinite chain of planar faces $F_n$ spanned by the vectors $v_n,v_{n+1}$, $n \in \mathbb Z$.

\smallskip

The $SL(3,\mathbb R)$-orbit of self-associated cones corresponding to $(a,b)$ can be parameterized by the initial condition $W$. }

\begin{figure}
\centering
\includegraphics[width = 14.92cm,height = 20.15cm]{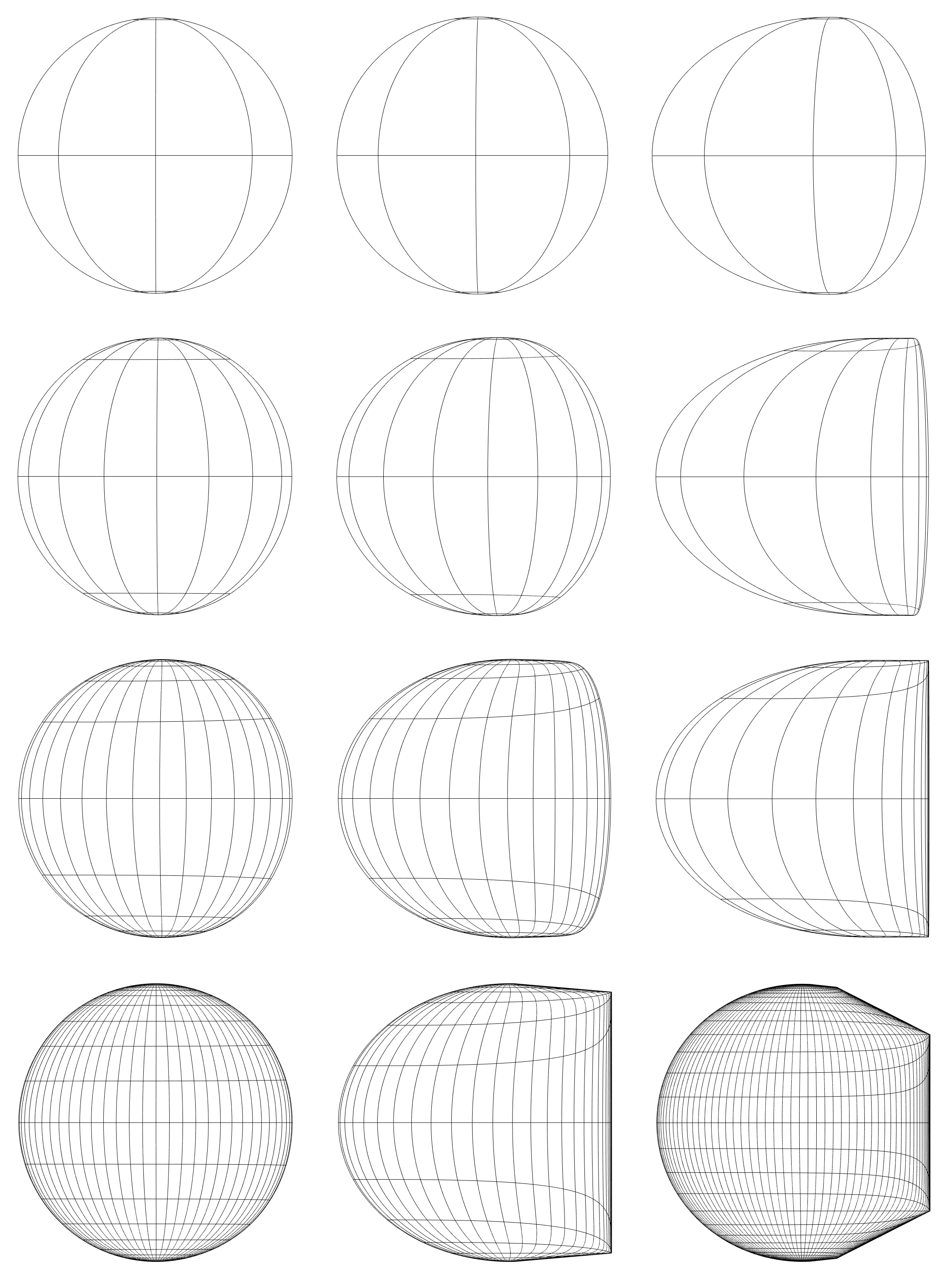}
\caption{Compact affine sections of self-associated cones of hyperbolic type with different parameter values $(a,b)$. The intervals $(a,b)$ are $(-3,2)$, $(-1,0)$, $(1,2)$ in the first row, $(-4,-2)$, $(-2,0)$, $(0,2)$ in the second row, $(-6,2)$, $(-4,0)$, $(-2,2)$ in the third row, and $(-12,-4)$, $(-6,2)$, $(-14,2)$ in the last row. The domain $M = (a,b) + i\mathbb R$ with uniform grid in cartesian coordinates is projected onto the interior of the section. The step size in the horizontal direction equals $\frac14$, in the vertical direction $\frac{\pi}{4}$. }
\label{fig_hyp}
\end{figure}

In Fig.~\ref{fig_hyp} we represent compact affine sections of self-associated cones of hyperbolic type for different parameter values $a,b$. As in the previous cases we project a uniformly spaced grid in the domain $M$ on the affine section.

\section{Open problems} \label{sec:conclusion}



A cone $K \subset \mathbb R^3$ is self-associated if whenever there is another cone $\tilde K \subset \mathbb R^3$ such that the affine spheres which are asymptotic to the boundaries $\partial K,\partial\tilde K$ are isometric, the cones $K,\tilde K$ are linearly isomorphic. This characterization of self-associated cones can be generalized to higher dimensions and suggests the following problem.

\medskip

\emph{Problem 1:} For $n > 3$, which cones $K \subset \mathbb R^n$ satisfy the following condition: whenever there exists another cone $\tilde K \subset \mathbb R^n$ such that the complete hyperbolic affine spheres with mean curvature $H = -1$ which are asymptotic to the boundaries $\partial K,\partial\tilde K$ are isometric, the cones $K,\tilde K$ are linearly isomorphic?

\medskip

The affine spheres which are asymptotic to the boundary of self-associated cones possess a continuous group of isometries which multiply the cubic differential representing the cubic form by unimodular complex constants different from 1, and hence do not preserve the cubic form. We may then ask whether this is possible in higher dimensions.

\medskip

\emph{Problem 2:} For $n > 3$, do there exist complete hyperbolic affine spheres $f: M \to \mathbb R^n$ with continuous groups of isometries which do not preserve the cubic form?

\medskip

Any associated family of $SL(3,\mathbb R)$-orbits of cones $K \subset \mathbb R^3$ admits a na\-tu\-ral action of the circle group $S^1$, multiplying the holomorphic function $U$ representing the cubic differential by unimodular complex constants. For a self-associated cone its associated family consists of a single orbit. On the other hand, for a cone which is merely linearly isomorphic to its dual cone by a unimodular isomorphism we have that the action of the element $e^{i\pi} \in S^1$ leaves the orbits in the associated family fixed. This is a weaker condition than being self-associated. We may then consider the following intermediate notions.

\medskip

\emph{Problem 3:} Let $k \geq 3$ be an integer. Find cones $K \subset \mathbb R^3$ (other than the self-associated cones) such that the action of the element $e^{2\pi i/k} \in S^1$ on its associated family of $SL(3,\mathbb R)$-orbits of cones leaves the orbits fixed.

\medskip

Affine spheres can be viewed as minimal Lagrangian manifolds in a certain para-K\"ahler space form \cite{Hildebrand11B}. The affine spheres corresponding to self-associated cones can hence be represented as minimal Lagrangian surfaces with a continuous symmetry group. Similar surfaces in the space $\mathbb {CP}^2$ have been considered in \cite{DorfmeisterMa16b},\cite{DorfmeisterMa16a} with loop group methods \cite{DorfmeisterEitner01}. Surfaces with rotational and with translational symmetries have been distinguished and their loop group decompositions have been computed. Loop group methods are applicable to the case of definite affine 2-spheres as well \cite{LiangJi10},\cite{LinWangWang17}.

\medskip

\emph{Problem 4:} Find the loop group decompositions of the affine spheres corresponding to the self-associated cones.

%
%
%
%
%

\section*{Acknowledgements}

The author would like to thank Prof.~Josef Dorfmeister for many inspiring discussions. In particular, he drew our attention to the fact that the reduction of Wang's equation is equivalent to the Painlev\'e III equation.

\appendix

\section{The degenerate Painlev\'e equation} \label{sec:PainleveIII}

In this section we consider the Painlev\'e III equation
\begin{equation} \label{PainleveIII}
tvv'' = t(v')^2 - vv' + v^3 - t
\end{equation}
with parameters $(\alpha,\beta,\gamma,\delta) = (1,0,0,-1)$, corresponding to the degenerate case $D_7$ in the classification of \cite{Sakai01}. The ensemble of its solutions $v_{s,c}(t)$ can be parameterized by two complex parameters $s,c$ which encode monodromy data of an associated linear system of ODEs \cite{ItsNovokshenov}. Equation \eqref{PainleveIII} has no solution expressible in terms of algebraic and classical transcendental functions except $v_{0,0}(t) = t^{1/3}$ \cite{OKSO06}. Equation \eqref{PainleveIII} has been studied in \cite{Kitaev87}. In particular, the asymptotics of the solutions in the neighbourhood of the singular points $t = 0$ and $t = +\infty$ have been computed. We summarize these results as follows.

{\prop \label{prop:Kitaev} The real solutions which are pole-free and positive in a neighbourhood of $t = 0$ on the positive real axis are parameterized by $(s,c) \in [-1,3] \times \mathbb R$, or equivalently $(\lambda,c) \in [0,1] \times \mathbb R$, where $\frac{s-1}{2} = \cos(\lambda\pi)$, and have the following asymptotics. For $s = 3$ or $\lambda = 0$
\[ v_{3,c}(t)t \sim \frac{2}{\left( \log t - 3\log 2 + \frac32\gamma + \frac34c \right)^2},
\]
for $0 < s < 3$ or $\lambda \in (0,\frac23)$
\[ v_{s,c}(t)t \sim \frac{2\lambda^2}{\sinh^2\left(\lambda\log t - 3\lambda\log 2 + \frac12\log\frac{\Gamma(1-\frac{\lambda}{2})\Gamma(1-\lambda)}{\Gamma(1+\frac{\lambda}{2})\Gamma(1+\lambda)} + \frac{3\lambda}{4}c\right)},
\]
for $s = 0$ or $\lambda = \frac23$
\[ v_{0,c}(t)t^{-1/3} \sim e^c \left(1 + \frac{9(e^c - e^{-2c})}{16}t^{4/3}\right),
\]
for $-1 < s < 0$ or $\lambda \in (\frac23,1)$
\[ v_{s,c}(t)t^{-1} \sim \frac{1}{2-2\lambda}\sinh\left(-(2-2\lambda)\log t + (6-6\lambda)\log 2 - \log\frac{-\Gamma(\frac{\lambda}{2})\Gamma(-1+\lambda)}{\Gamma(1-\frac{\lambda}{2})\Gamma(1-\lambda)} + 3(1-\lambda)c\right),
\]
for $s = -1$ or $\lambda = 1$
\[ v_{-1,c}(t)t^{-1} \sim -\log t + 2\log 2 - \frac32\gamma + \frac32c,
\]
where $\gamma$ is the Euler-Mascheroni constant\footnote{In the original paper \cite{Kitaev87} a factor of 2 is missing in the last expression.}.

The real solutions which are pole-free and positive in a neighbourhood of $t = +\infty$ on the positive real axis are parameterized by $s \in \mathbb R$, the second parameter $c$ being zero, and have the asymptotics
\[ v_{s,0}(t)t^{-1/3} - 1 \sim 3^{-1/4}\pi^{-1/2}st^{-1/3}e^{-3\sqrt{3}t^{2/3}/2}.
\] }

The solutions $v_{s,0}$ with $s \in [-1,3]$ appear in both lists, and their asymptotics both at $t = 0$ and $t = +\infty$ are known.

For studying positive solutions $v(t)$ of \eqref{PainleveIII} it is convenient to make the substitution $\tau = \log t$, $g = \log v - \frac13\log t$. Then \eqref{PainleveIII} is equivalent to the second-order ODE
\begin{equation} \label{g_tau_eq}
\frac{d^2g}{d\tau^2} = e^{4\tau/3}(e^g - e^{-2g})
\end{equation}
on the function $g(\tau)$. By virtue of the strict monotonicity of the function $e^g - e^{-2g}$ we have $g'' > \tilde g''$ for two solutions $g,\tilde g$ whenever $g > \tilde g$.

We have the following monotonicity results.

{\lemma \label{lem:s0monotone} Let $t_a \geq 0$ and $s < s'$ be real numbers such that the solutions $v_{s,0},v_{s',0}$ of \eqref{PainleveIII} are positive on the interval $(t_a,+\infty)$. Let $g_s,g_{s'}$ be the corresponding solutions of \eqref{g_tau_eq}. Then $g_{s'} > g_s$, $g'_{s'} < g'_s$ on $(\log t_a,+\infty)$ and $v_{s,0} < v_{s',0}$ on $(t_a,+\infty)$. 

In particular, if the solution $v_{s,0}$ is positive on $\mathbb R_+$, then for $s' > s$ the solution $v_{s',0}$ is positive either on $\mathbb R_+$ or up to the right-most pole, and for $s' < s$ the solution $v_{s',0}$ must cease to be positive before it reaches its right-most pole. }

\begin{proof}
The asymptotics of $v_{s,0}$ for $t \to +\infty$ yields
\[ g_s(\tau) \sim 3^{-1/4}\pi^{-1/2}se^{-\tau/3-3\sqrt{3}\exp(2\tau/3)/2}
\]
as $\tau \to +\infty$. Define the function $\delta = g_{s'} - g_s$ on $(\log t_a,+\infty)$. We get
\[ \delta(\tau) \sim 3^{-1/4}\pi^{-1/2}(s'-s)e^{-\tau/3-3\sqrt{3}\exp(2\tau/3)/2}
\]
as $\tau \to +\infty$. Moreover, $\delta'' > 0$ 
for all $\tau$ such that $\delta(\tau) > 0$.

Now we have $\delta(\tau) > 0$ for large enough $\tau$, and $\lim_{\tau\to+\infty}\delta(\tau) = 0$. Hence there exists a sequence $\tau_k \to +\infty$ such that $\delta(\tau_k) > 0$, $\delta'(\tau_k) \leq 0$ for all $k$. But then $\delta''(\tau) > 0$, $\delta'(\tau) < 0$, $\delta(\tau) > 0$ for all $\tau < \tau_k$, and hence for all $\tau \in (\log t_a,+\infty)$. The first claim of the lemma now easily follows.

Let us prove the second claim. For $s' > s$ we have that $v_{s',0}(t) > v_{s,0}(t) > 0$ as long as $v_{s',0} > 0$ on $(t,+\infty)$. Therefore the only way for $v_{s',0}$ to cease to be positive is to cross a pole. For $s' < s$ we have $v_{s',0}(t) < v_{s,0}(t)$ as long as $v_{s',0} > 0$ on $(t,+\infty)$. Hence $v_{s',0}$ must first hit the zero value before it can reach a pole.
\end{proof}

{\lemma \label{lem:c_monotone} Let $s \in [0,3]$, $t_b > 0$, and $c < c'$ be real numbers such that the solutions $v_{s,c},v_{s,c'}$ of \eqref{PainleveIII} are positive on the interval $(0,t_b)$. Let $g_c,g_{c'}$ be the corresponding solutions of \eqref{g_tau_eq}. Then $g_{c'} > g_c$, $g'_{c'} > g'_c$ on $(-\infty,\log t_b)$ and $v_{s,c} < v_{s,c'}$ on $(0,t_b)$. 

In particular, if the solution $v_{s,0}$ is positive on $\mathbb R_+$, then for $c > 0$ the solution $v_{s,c}$ is positive up to the left-most pole, and for $c < 0$ the solution $v_{s,c}$ cannot be positive up to its left-most pole. }

\begin{proof}
The asymptotics of $v_{s,c}$ for $t \to 0$ from Proposition \ref{prop:Kitaev} yields
\[ g_c(\tau) \sim \left\{ \begin{array}{rcl} -\frac43\tau + \log 2 - 2\log(-\tau + 3\log 2 - \frac32\gamma - \frac34c),& \quad & s = 3; \\ 
(-\frac43+2\lambda)\tau + 2\log\lambda + 3(1 - 2\lambda)\log 2 + \log\frac{\Gamma(1-\frac{\lambda}{2})\Gamma(1-\lambda)}{\Gamma(1+\frac{\lambda}{2})\Gamma(1+\lambda)} + \frac{3\lambda}{2}c,&& s \in (0,3); \\
c + \frac{9(e^c - e^{-2c})}{16}e^{4\tau/3},&& s = 0 \end{array} \right.
\]
as $\tau \to -\infty$. Here $\frac{s-1}{2} = \cos(\pi\lambda)$. Define the function $\delta = g_{c'} - g_c$ on $(-\infty,\log t_b)$. We get
\[ \delta(\tau) \sim \left\{ \begin{array}{rcl} \frac{3(c'-c)}{-2\tau},& \quad & s = 3; \\
\frac{3\lambda(c'-c)}{2},&& s \in [0,3) \end{array} \right.
\]
as $\tau \to -\infty$. Recall that $\delta'' > 0$ for all $\tau$ such that $\delta > 0$.

It follows that $\delta > 0$, $\delta'' > 0$ for large enough $|\tau|$, and $\lim_{\tau\to-\infty}\delta(\tau) < \infty$. But then also $\delta' > 0$ for large enough $|\tau|$. Therefore $\delta$ cannot become negative on the whole interval $(-\infty,\log t_b)$, and the first claim of the lemma follows.

In order to prove the second claim, compare the solution $v_{s,c}$ to the solution $v_{s,0}$. For $c > 0$ we have that $v_{s,c}(t) > v_{s,0}(t) > 0$ as long as $v_{s,c} > 0$ on $(0,t)$. Therefore the only way for $v_{s,c}$ to cease to be positive is to cross a pole. For $c < 0$ we have $v_{s,c}(t) < v_{s,0}(t)$ as long as $v_{s,c} > 0$ on $(0,t)$. Therefore the only way for $v_{s,c}$ to cease to be positive is to become zero. 
\end{proof}

By the Painlev\'e property the solutions of \eqref{PainleveIII} can be extended to meromorphic functions on the universal cover of $\mathbb C \setminus \{0\}$. If $t = t_0 \not= 0$ is a singularity of the solution, then inserting its Laurent expansion around $t_0$ into \eqref{PainleveIII} easily yields
\begin{equation} \label{PIIIpole_expansion}
v(t) = \frac{2t_0}{(t-t_0)^2} + \alpha - \frac{\alpha}{t_0}(t-t_0) + \frac{3\alpha^2t_0 + 9\alpha}{10t_0^2}(t-t_0)^2 - \frac{3\alpha^2t_0 + 4\alpha}{5t_0^3}(t-t_0)^3 + \dots,
\end{equation}
for some $\alpha \in \mathbb R$, i.e., $t_0$ is a double pole of the solution. 

The family $(v_{\alpha})_{\alpha \in \mathbb R}$ of solutions having a double pole at a fixed $t_0$ satisfies a similar monotonicity result.

{\lemma \label{lem:alpha_monotone} Let $t_0 > 0$, $\alpha < \alpha'$ be real numbers, and let $v_{\alpha},v_{\alpha'}$ be the solutions of \eqref{PainleveIII} with the corresponding expansions \eqref{PIIIpole_expansion}. Suppose that $v_{\alpha},v_{\alpha'}$ are positive on the interval $(T,t_0)$ for some $T \in [-\infty,t_0)$, and let $g_{\alpha},g_{\alpha'}$ be the corresponding solutions of \eqref{g_tau_eq} on $(\log T,\tau_0)$, where $\tau_0 = \log t_0$. Then we have $g_{\alpha} < g_{\alpha'}$, $g'_{\alpha} > g'_{\alpha'}$, $g''_{\alpha} < g''_{\alpha'}$ on $(\log T,\tau_0)$, and $v_{\alpha} < v_{\alpha'}$ on $(T,t_0)$. }

\begin{proof}
Set $\delta = g_{\alpha'} - g_{\alpha}$. Since the ratio $\frac{v_{\alpha'}}{v_{\alpha}}$ is analytic and non-zero in the neighbourhood of $t_0$, its logarithm is also analytic and by virtue of \eqref{PIIIpole_expansion} has the expansion
\[ \delta(\tau_0+d) = \frac{\alpha'-\alpha}{2t_0}d^2 - \frac{\alpha' - \alpha}{2t_0^2}d^3 + \frac{t_0((\alpha')^2-\alpha^2) + 18(\alpha'-\alpha)}{40t_0^3}d^4 - \frac{t_0((\alpha')^2-\alpha^2) + 8(\alpha'-\alpha)}{20t_0^4}d^5 + \dots.
\]
It follows that for $\tau < \tau_0$ close enough to $\tau_0$ we have $\delta(\tau) > 0$ and $\delta'(\tau) < 0$. But $\delta'' > 0$ as long as $\delta > 0$, hence these relations are valid on the whole interval $(\tau_0,\log T)$. The claim of the lemma now easily follows.
\end{proof}

{\corollary \label{cor:s_unique} Let $s,s' \in (-1,3)$ and $c,c' \in \mathbb R$ be numbers, and suppose that the corresponding solutions $v_{s,c},v_{s',c'}$ of \eqref{PainleveIII} are positive on $(0,t_0)$ and have a double pole at $t_0$, such that the corresponding constants $\alpha,\alpha'$ in their expansions \eqref{PIIIpole_expansion} satisfy $\alpha < \alpha'$. Then $s < s'$. }

\begin{proof}
Let $g_{s,c},g_{s',c'}$ be the corresponding solutions of \eqref{g_tau_eq}, and set $\delta = g_{s',c'} - g_{s,c}$. By Lemma \ref{lem:alpha_monotone} we have $\delta' < 0$, $\delta'' > 0$ on $(-\infty,\log\tau_0)$. From Proposition \ref{prop:Kitaev} we get after some calculations that the asymptotics of $g_{s,c},g_{s',c'}$ in the limit $\tau \to -\infty$ are given by
\[ g_{s,c} \sim (2\lambda - \frac43)\tau + C(s,c), \quad g_{s',c'} \sim (2\lambda' - \frac43)\tau + C(s',c'),
\]
where $\frac{s-1}{2} = \cos(\pi \lambda)$, $\frac{s'-1}{2} = \cos(\pi \lambda')$, and $C(s,c),C(s',c')$ are constants. It follows that $\delta \sim 2(\lambda' - \lambda)\tau + const$ as $\tau \to -\infty$. But $\delta$ is a strictly decreasing convex function, and therefore we must have $\lambda' < \lambda$. This implies $s' > s$.
\end{proof}

Finally, we shall estimate the error in the asymptotics of the solution $v_{3,c}$ provided in \cite{Kitaev87}. For given $c \in \mathbb R$, set for brevity $c' = 3\log 2 - \frac32\gamma - \frac34c$.

{\lemma \label{lem:asymp_v3c} The asymptotics of the solution $v_{3,c}$ for $t \to 0$ is given by $\frac{2t^{-1}}{( - \log t + c' )^2} - \frac{t^3(-\log t+c'+1)^2}{32}$. The asymptotics of its derivative is given by $\frac{dv_{3,c}}{dt} \sim -\frac{2t^{-2}}{( -\log t + c' )^2} + \frac{4t^{-2}}{( -\log t + c' )^3} - \frac{3t^2(-\log t + c' + \frac23)^2}{32}$. }

\begin{proof}
By Proposition \ref{prop:Kitaev} the solution $g(\tau)$ of \eqref{g_tau_eq} which corresponds to $v_{3,c}$ has asymptotics
\[ g(\tau) \sim \log2 - \frac43\tau - 2\log(-\tau+c') =: \tilde g(\tau).
\]
Note that the function $\tilde g$ satisfies the ODE $\tilde g'' = e^{4\tau/3}e^{\tilde g}$, and hence $\delta = \tilde g - g$ satisfies
\begin{equation} \label{delta_ODE}
\delta'' = e^{4\tau/3}(e^{\tilde g} - e^g + e^{-2g}) = 2\left(-\tau+c'\right)^{-2}(1 - e^{-\delta}) + \mu(\tau),
\end{equation}
where $\mu(\tau) = e^{4\tau/3}e^{-2g} \sim \frac14e^{4\tau}(-\tau+c')^4$.

Define the functions $\underline{\mu}(\tau) = \int_{-\infty}^{\tau}\int_{-\infty}^s\mu(r)\,dr\,ds$ and $\overline{\mu}(\tau) = \frac{1}{-\tau+c'}\int_{-\infty}^{\tau}(s-c')^2\int_{-\infty}^s\frac{\mu(r)}{-r+c'}\,dr\,ds$. Note that $\overline{\mu}'' = \frac{2\overline{\mu}}{(\tau-c')^2} + \mu > \frac{2}{(\tau-c')^2}(1 - e^{-\overline{\mu}}) + \mu$. Note that both $\underline{\mu}$ and $\overline{\mu}$ have the same asymptotics $\frac{e^{4\tau}(-\tau+c'+\frac12)^4}{64}$ as $\tau \to -\infty$. We shall now show that $\underline{\mu} < \delta < \overline{\mu}$.

Let $\tau_0 < c'$ be such that $g$ and hence also $\mu$ and $\delta$ are defined on $(-\infty,\tau_0]$. Consider the integral operator ${\cal I}$ taking a smooth function $f: (-\infty,\tau_0] \to \mathbb R$ to $[{\cal I}f](\tau) = \int_{-\infty}^{\tau}\int_{-\infty}^s \frac{2}{(-r+c')^2}(1 - e^{-f(r)}) + \mu(r)\,dr\,ds$. If $f$ is such that $0 \leq f(\tau) \leq \overline{\mu}(\tau)$ for all $\tau \in (-\infty,\tau_0]$, then
\[ 0 \leq [{\cal I}f](\tau) \leq \int_{-\infty}^{\tau}\int_{-\infty}^s \frac{2}{(-r+c')^2}(1 - e^{-\overline{\mu}(r)}) + \mu(r)\,dr\,ds < \int_{-\infty}^{\tau}\int_{-\infty}^s \overline{\mu}''(r)\,dr\,ds = \overline{\mu}(\tau), \quad \tau \leq \tau_0.
\]
Hence ${\cal I}$ is well-defined for every such function $f$. Clearly if $f < \tilde f$ everywhere on $(-\infty,\tau_0]$, then also ${\cal I}f < {\cal I}\tilde f$ everywhere on this interval.

Define now $\delta_0 \equiv 0$ and recursively $\delta_{k+1} = {\cal I}\delta_k$, $k \geq 0$. Then $0 = \delta_0 < \delta_1 = \underline{\mu}$, and hence the sequence $\delta_k(\tau)$ is strictly increasing for every $\tau \leq \tau_0$. On the other hand, this sequence is upper bounded by $\overline{\mu}(\tau)$. Hence $\delta_k$ converges point-wise to some limit function $\delta^*$. This function is a fixed point of the operator ${\cal I}$ and hence smooth and a solution of ODE \eqref{delta_ODE}. Thus it coincides with the sought function $\delta$. By construction we obtain the desired bounds $\underline{\mu} < \delta < \overline{\mu}$ and therefore the asymptotics $\delta \sim \frac{e^{4\tau}(-\tau+c'+\frac12)^4}{64}$.

Switching back to the solution $v(t)$ of \eqref{PainleveIII} we get
\[ v_{3,c}(t) = t^{1/3}e^{\tilde g - \delta} \sim \frac{2t^{-1}}{( - \log t + c' )^2}\left(1 - \frac{t^4(-\log t+c'+\frac12)^4}{64} \right),
\]
which yields the desired error term for $v_{3,c}$.

Further, by \eqref{delta_ODE} we have $\delta'' \sim \frac{e^{4\tau}(-\tau+c')^4}{4}$ and hence $\delta'(\tau) = \int_{-\infty}^{\tau}\delta''(s)\,ds \sim \frac{e^{4\tau}(-\tau+c'+\frac14)^4}{16}$. It follows that
\[ g'(\tau) = \tilde g' - \delta' \sim - \frac43 + \frac{2}{-\tau+c'} - \frac{e^{4\tau}(-\tau+c'+\frac14)^4}{16}.
\]
Therefore by virtue of $\frac{d\tau}{dt} = t^{-1}$ we have
\begin{align*}
\frac{dv_{3,c}}{dt} &= t^{-2/3}e^{\tilde g - \delta}\left(g' + \frac13\right) \\
&\sim \frac{2t^{-2}}{(-\log t+c')^2}\left( 1 - \frac{t^4(-\log t+c'+\frac12)^4}{64} \right)\left(- 1 + \frac{2}{-\log t+c'} - \frac{t^4(-\log t+c'+\frac14)^4}{16} \right),
\end{align*}
which after a little calculus yields the desired error term.
\end{proof}

%

\bibliographystyle{plain}
\bibliography{convexity,geometry,affine_geometry,diff_eq,misc}

\end{document}